\documentclass{article} 
\usepackage{iclr2023_workshop,times}
\usepackage{graphicx}

\author{Nick McGreivy \\
Department of Astrophysical Sciences\\
Princeton University\\
\texttt{mcgreivy@princeton.edu}
\And
Ammar Hakim \\
Princeton Plasma Physics Laboratory \\
Princeton, NJ \\
\texttt{ahakim@pppl.gov}
}

\iclrfinalcopy

\usepackage{amssymb}
\usepackage{bm}
\usepackage{amsmath}
\usepackage{mathtools}
\usepackage[backref=page]{hyperref}
\hypersetup{
           breaklinks=true,   
}
\usepackage{cleveref}
\usepackage{nicefrac}
\usepackage{subcaption}

\title{{Invariant Preservation in Machine Learned PDE Solvers via Error Correction}}

\begin{document}

\maketitle

\begin{abstract}
\looseness=-1 Machine learned partial differential equation (PDE) solvers trade the reliability of standard numerical methods for potential gains in accuracy and/or speed.
The only way for a solver to guarantee that it outputs the exact solution is to use a convergent method in the limit that the grid spacing $\Delta x$ and timestep $\Delta t$ approach zero.
Machine learned solvers, which learn to update the solution at large $\Delta x$ and/or $\Delta t$, can never guarantee perfect accuracy.
Some amount of error is inevitable, so the question becomes: how do we constrain machine learned solvers to give us the sorts of errors that we are willing to tolerate?
In this paper, we design more reliable machine learned PDE solvers by preserving discrete analogues of the continuous invariants of the underlying PDE.
Examples of such invariants include conservation of mass, conservation of energy, the second law of thermodynamics, and/or non-negative density.
Our key insight is simple: to preserve invariants, at each timestep apply an error-correcting algorithm to the update rule.
Though this strategy is different from how standard solvers preserve invariants, it is necessary to retain the flexibility that allows machine learned solvers to be accurate at large $\Delta x$ and/or $\Delta t$.
This strategy can be applied to any autoregressive solver for any time-dependent PDE in arbitrary geometries with arbitrary boundary conditions.
Although this strategy is very general, the specific error-correcting algorithms need to be tailored to the invariants of the underlying equations as well as to the solution representation and time-stepping scheme of the solver. The error-correcting algorithms we introduce have two key properties. First, by preserving the right invariants they guarantee numerical stability. Second, in closed or periodic systems they do so without degrading the accuracy of an already-accurate solver.

\end{abstract}

\section{Introduction}\label{sec:intro}

Scientists and engineers are interested in solving partial differential equations (PDEs). Many PDEs cannot be solved analytically, and must be approximated using discrete numerical algorithms. We refer to these algorithms as `PDE solvers.' The fundamental challenge for PDE solvers is to balance between two competing objectives: first, to find an accurate approximation to the solution of the equation, and second, to do so with as few computational resources as possible.

Decades of research into discrete numerical algorithms have resulted in reliable solvers for most PDEs of interest. For time-dependent PDEs, these so-called `standard numerical methods' use hand-crafted rules to update the solution at each timestep. 
Successful hand-crafted update rules have two key properties. First, the property of \textit{convergence}. Convergent methods converge to the exact solution in the limit that the grid spacing $\Delta x$ and the timestep $\Delta t$ approach zero \cite{lax1956survey}. Second, the property of \textit{invariant preservation}. Time-dependent PDEs often have one or more invariants. Examples of such invariants include conservation of mass \cite{eymard2000finite}, conservation of energy \cite{hakim2019discontinuous}, non-decreasing entropy \cite{merriam1989entropy}, and/or non-negative density \cite{perthame1996positivity}. Invariant preserving PDE solvers satisfy discrete analogues of these continuous invariants when $\Delta x$ and $\Delta t$ are positive. As a result, they are numerically stable and do not violate qualitatively important properties.

In recent years, scientists and engineers have attempted to use machine learning (ML) to design new and better PDE solvers  \cite{accelerating_eulerian_fluid_sim,data_driven_discretizations,2d_learned_advection,embedding_hard_constraints,learning_neural_pde_solvers,rose_yu_paper,les_closure,solver-in-the-loop, mishra2018machine}.
The goal of these machine learned PDE solvers is as follows. Suppose there is a PDE we would like to find an approximate solution to for many different initial or boundary conditions (ICs or BCs) or on a very large domain. We first generate training data, usually from a highly accurate standard solver. Next, we design and train a learnable update rule; often this involves the use of neural networks. If the learned update rule is faster than the best-performing standard solvers with equal accuracy, it can then be used to amortize the initial cost of training over many different ICs or BCs or a much larger domain, by finding sufficiently accurate solutions at reduced computational cost \cite{mishra2018machine,xue2020amortized}. 

For time-dependent PDEs, the main strategy for making solvers faster with ML is to solve the equations at coarser resolution, i.e., with increased $\Delta x$ and/or $\Delta t$ relative to standard methods \cite{vinuesa2022enhancing}.\footnote{There are other possibilities for accelerating the solution of PDEs with ML. For stiff PDEs, it may be possible to replace a solver that requires implicit (slow) timestepping with an explicit machine learned surrogate \cite{wang2023long,dileoni2021deeponet}. For other PDEs, it may be possible to accelerate a solver by learning a surrogate for some particularly costly operator at fixed grid resolution. One can, for example, use ML to accelerate a Poisson solve \cite{cheng2021using,pmlr-v97-greenfeld19a} or learn to approximate a costly plasma collision operator \cite{miller2021encoder,dener2020training,holloway2021acceleration}.}
To do so, machine learned solvers use a radically different approach from standard PDE solvers. Instead of designing hand-crafted update rules that converge as $\Delta x \rightarrow 0$ and $\Delta t \rightarrow 0$, machine learned solvers learn an update rule from data that is accurate at some large value(s) of $\Delta x$ and/or $\Delta t$.\footnote{Note that for time-\textit{independent} PDEs, it is possible to guarantee convergence by interfacing ML with a standard iterative solver. See, for example, \cite{illarramendi2022performance}.} Ideally, this update rule would be faster than standard methods at a given level of accuracy.
Machine learned PDE solvers were found, for certain problems, to have successfully achieved high accuracy at low computational cost \cite{ml_accelerated_cfd, deepmind_turbulence_sims, fourier_neural_operator, multigrid, multigrid_graph_nn}.
A more recent study \cite{dresdner2022learning}, which makes a comparison between a very efficient standard numerical method and the best ML-based solvers, tacitly implies that earlier studies compared to weak baselines and thus casts doubt on some purportedly impressive claims about improved performance from ML-based solvers. Nevertheless, there is still potential for ML-based solvers to accelerate the solution of PDEs, and this remains an active area of research.

Because the only way for a numerical method to guarantee perfect accuracy (i.e., to output the exact solution) is (except in trivial cases) to use a convergent method in the limit $\Delta x \rightarrow 0$ and $\Delta t \rightarrow 0$, there is no way to guarantee that machine learned solvers output the exact solution while using large $\Delta x$ and/or $\Delta t$.
Some amount of error is inevitable, so the question becomes: how can we constrain machine learned solvers to give us the sorts of errors that we are willing to tolerate? In other words, how can we build more reliable machine learned PDE solvers?

One approach to building more reliable machine learned PDE solvers is to change the ML model and training procedure. This approach is quite natural to students of ML. Improving the model \cite{message_passing_neural_pde_solvers}, increasing the dataset size \cite{brandstetter2022lie}, changing the loss function \cite{solver-in-the-loop}, adding random noise to the training procedure \cite{deepmind_turbulence_sims,sanchez2020learning}, and adding regularization \cite{Kaptanoglu_2021} are all examples of this approach. To some extent, these techniques have been successful at improving robustness and numerical stability \cite{klimesch2022simulating}. However, none of these ML-based techniques are capable of \textit{guaranteeing} numerical stability. While these solvers may give reliable results for some inputs, on other inputs the solution might blow up or be nonsensical.

The purpose of this paper is to propose and demonstrate a different and mutually compatible approach to building more reliable machine learned PDE solvers: \textit{preserving discrete invariants}. This approach is quite natural to students of computational physics, as so much of the theory and development of standard numerical methods is related to ensuring those methods preserve the right invariants. This approach can be used with any solver that uses an update rule (including standard solvers) and is otherwise agnostic to the details of the solver.

Why do we want our machine learned solvers to preserve invariant quantities? A simple but incomplete explanation is that ML models which use physical knowledge as an inductive bias tend to outperform models that don't \cite{thuerey_guaranteed_momentum,ling2016reynolds,brandstetter2022clifford}. Invariant preservation is, when done correctly, free lunch. We know that for a given PDE our solution should preserve certain invariants, so by enforcing those invariants at each timestep we improve the solution. A second reason has to do with \textit{numerical stability}. By preserving the right combination of invariants, we can design machine learned PDE solvers which are numerically stable by construction. These solvers, like well-designed standard solvers, are guaranteed not to blow up as $t \rightarrow \infty$. A third reason has to do with \textit{trust}.
People are unlikely to use solvers they do not trust.
People are more likely to trust a numerical method if it preserves the correct set of invariants.

In standard numerical methods, most of the theory of invariant preservation has been developed for hyperbolic PDEs in conservation form
\begin{equation}\label{eq:hyperbolic}
    \frac{\partial \bm u}{\partial t} + \bm \nabla \cdot \bm F(\bm u) = 0
\end{equation}
where $\bm u(\bm x, t) \in \mathbb{R}^m$, $\bm x \in \Omega$, the domain $\Omega \subset \mathbb{R}^d$, $t \in \mathbb{R}^+$, and $\bm F(\bm u) \in C^1(\mathbb{R}^{m \times d})$. This is because hyberbolic PDEs are conservation laws which maintain a variety of conserved quantities and other physical invariants. Because hyperbolic PDEs are so important to research and production applications, because their theory is so well-developed, and because they make a useful playground for studying invariants in time-dependent PDEs, this paper is focused on preserving invariants in machine learned solvers for hyberbolic PDEs. However, as we discuss in \cref{sec:mlpdesolvers}, the strategy and techniques we introduce can be applied to other invariant-preserving PDEs as well as PDEs with non-invariant terms.

The key insight of this paper is simple: to preserve discrete invariants in machine learned PDE solvers, at each timestep apply an error-correcting algorithm to the update rule. Let us now sketch how this works. Suppose that we represent the continuous solution $\bm u(\bm x, t)\in \mathbb{R}^m$ to \cref{eq:hyperbolic} with a discrete solution $\hat{\bm u}(\bm x, t) \in \mathbb{R}^m$ which is a linear sum of $N$ basis functions $\bm \phi_k(\bm x) \in \mathbb{R}^m$ and coefficients $\bm c_k(t) \in \mathbb{R}^m$ for $k \in [1, \dots, N]$, such that $\hat{u}_j(\bm x, t) = \sum_{k=1}^N c_{jk}(t) \phi_{jk}(\bm x)$ for $j \in [1, \dots, m]$. 
Suppose also that the update rule predicts that $\hat{\bm u}$ will change at a rate $\frac{\partial \hat{\bm u}}{\partial t}$ at time $t$ (or an amount $\Delta \hat{\bm u}$ in time $\Delta t$). Suppose also that, due to a priori knowledge about the underlying equations, we would like the solution to satisfy $L$ discrete invariants $\mathcal{I}_\ell(\hat{\bm u}, \frac{\partial \hat{\bm u}}{\partial t})$ for $\ell = [1, \dots, L]$ satisfying either equalities ($\mathcal{I}_\ell= 0$) or inequalities ($\mathcal{I}_\ell \ge 0$). 
If $\frac{\partial \hat{\bm u}}{\partial t}$ (or $\Delta \hat{\bm{u}})$ does not already satisfy these discrete invariants, then we use an error correcting algorithm to modify $\frac{\partial \hat{\bm u}}{\partial t}$ (or $\Delta \hat{\bm{u}})$ to ensure that each of the invariants is satisfied. We repeat this process each timestep. This process is meant to used at inference, though it could also be used while training.

Students of computational physics will notice that this strategy is very different from how invariants are preserved in standard solvers. As we discuss in \cref{sec:standardsolvers}, standard solvers typically conserve linear invariants by using some type of finite volume (FV) method where the rate of change of the discrete solution in each grid cell is equal to the flux through cell boundaries. Standard solvers typically conserve non-linear invariants by putting locally-derived constraints on the flux \cite{mishra2019numerical,hakim2019discontinuous}. Examples of such constraints include upwinding, centered fluxes, or limiters. Yet the error-correcting algorithms we introduce in \cref{sec:scalarhyperbolic,sec:systemshyperbolic} involve using \textit{global} constraints to constrain either the flux through cell boundaries or the update $\frac{\partial \hat{\bm u}}{\partial t}$ (or $\Delta \hat{\bm u}$). In these algorithms, instead of the time-derivative in a cell depending only on its nearest neighbors, the time-derivative depends on the value of the solution over the entire domain. 

One might ask why machine learned solvers should use global, rather than local, constraints to preserve invariants? It is \textit{possible} for a machine learned solver to preserve invariants by putting local constraints on the flux. The problem with doing so, as we show in \cref{sec:why_standard_dont_work}, is that for some invariants standard approaches can be too constraining -- in particular, they introduce too much numerical diffusion -- to give accurate results at large $\Delta x$, and are too constrained by the Courant-Friedrichs-Lewy (CFL) condition to work at large $\Delta t$. Global constraints, in contrast, allow machine learned solvers the flexibility to make accurate predictions at large $\Delta x$ and/or $\Delta t$, while ensuring that whatever errors are made at least do not violate invariants. 

The error-correcting algorithms we introduce have two key properties. First, by preserving the right invariants these algorithms guarantee numerical stability. Numerical stability is about whether the solution blows up; formally, this can be defined as $||\hat{\bm{u}}(\bm x, T)|| \le C_T ||\hat{\bm{u}}(\bm x, 0)||$ for time $T>0$ for some norm and some constant $C_T$ which can depend on $T$ but not on $\Delta x$ or $\Delta t$. For invariant-preserving PDEs, the inequality $||\hat{\bm{u}}(\bm x, T)|| \le ||\hat{\bm{u}}(\bm x, 0)||$ is often a better definition of stability \cite{dale_durran}. For scalar hyperbolic PDEs and for many systems of hyperbolic PDEs, it can be shown that preserving discrete analogues of a subset of the continuous invariants is sufficient to guarantee stability \cite{mishra2019numerical,merriam1989entropy,juno2018discontinuous}.

The second key property is that, as we show in \cref{sec:verification}, in closed or periodic systems this strategy preserves invariants without degrading the accuracy of an already-accurate solver. In closed or periodic systems, we know a priori the discrete invariants ($\mathcal{I}_\ell = 0$ or $\mathcal{I}_\ell \ge 0$). If the solver doesn't satisfy those invariants at each timestep, then we know that our solver has made an error. In open systems, we may not know the invariants due to uncertainty in the fluxes through the boundary. We can estimate these fluxes, but this estimate can be inaccurate. Thus, as we demonstrate in \cref{sec:verification_compressible_euler}, in open systems the error-correction process can degrade accuracy if the rate of change of the invariants cannot be estimated accurately.



Be careful not to confuse the property of invariance or equivariance under a transformation with the preservation of an invariant. Noether's theorem says that for certain physical systems, a continuous symmetry (invariance) leads to the existence of a conserved (invariant) quantity. However, for discrete systems no such theorem exists. While equivariant neural networks may be used to enforce discrete symmetries and thereby improve the generalization capabilities of machine learned PDE solvers \cite{wang2020incorporating}, whether or not a solver is equivariant is unrelated to whether or not that solver preserves an invariant.

The rest of this paper is structured as follows. 
In \cref{sec:mlpdesolvers}, we define what a machine learned PDE solver is and for what types of PDEs this strategy can be applied. 
In \cref{sec:standardsolvers}, we review the theory of invariant preservation in scalar hyperbolic PDEs and in systems of hyperbolic PDEs.
In \cref{sec:scalarhyperbolic}, we design invariant-preserving error-correcting algorithms for a variety of popular solver types for scalar hyperbolic PDEs.
In \cref{sec:systemshyperbolic}, we design invariant-preserving error-correcting algorithms for systems of hyperbolic PDEs, focusing on FV-like solvers. In \cref{sec:verification}, we turn our attention to machine learned solvers. We computationally verify that, at least in periodic systems, these error-correcting invariant-preserving algorithms do not degrade the accuracy of an already-accurate machine learned solver.   
Code and instructions for reproducing the figures in \cref{sec:scalarhyperbolic,sec:systemshyperbolic,sec:verification} can be found at \href{https://github.com/nickmcgreivy/InvariantPreservingMLSolvers}{https://github.com/nickmcgreivy/InvariantPreservingMLSolvers}.
\Cref{sec:relatedwork} is related work. We conclude in \cref{sec:limitations} by discussing limitations and trade-offs associated with invariant preservation.

\section{Machine Learned Partial Differential Equation Solvers}\label{sec:mlpdesolvers}

We are interested in finding numerical solutions to PDEs of the form
\begin{equation}\label{eq:generalpde}
    \frac{\partial \bm u}{\partial t} + \mathcal{N}[\bm u] = 0
\end{equation}
where the solution $\bm u(\bm x, t) \in \mathbb{R}^m$, $\bm x \in \Omega$, $\Omega \subset \mathbb{R}^d$, $t \in \mathbb{R}^+$, and $\mathcal{N} : \mathbb{R}^m \rightarrow \mathbb{R}^m$ is an operator which preserves some invariant quantities across time. Examples of PDEs in the form of \cref{eq:generalpde} include the parabolic diffusion equation
\begin{equation}
    \frac{\partial \bm u}{\partial t} = D{\nabla}^2 \bm u,
\end{equation}
scalar hyperbolic PDEs in the form
\begin{equation}
    \frac{\partial u}{\partial t} + \bm \nabla \cdot \bm f(u) = 0,
\end{equation}
Hamiltonian systems evolving under the influence of a Hamiltonian $H$ and Poisson bracket $\{\cdot,\cdot\}$ 
\begin{equation}\label{eq:hamiltoniansystem}
    \frac{\partial f}{\partial t} + \{ f, H\} = 0,
\end{equation}
systems of hyperbolic PDEs
\begin{equation}
    \frac{\partial \bm u}{\partial t} + \bm \nabla \cdot \bm F(\bm u) = 0,
\end{equation}
and the Boltzmann equation of kinetic physics which both advects particles in phase space and evolves them according to a mass, energy, and momentum-conserving collision operator
\begin{equation}
    \frac{\partial f}{\partial t} + \bm v \cdot \frac{\partial f}{\partial \bm x} + \bm a \cdot \frac{\partial f}{\partial \bm v} = \bigg(\frac{\partial f}{\partial t}\bigg)_{\textnormal{coll}}.
\end{equation}
To numerically approximate the solution to an equation in the form \cref{eq:generalpde}, we begin by representing the continuous solution $\bm u(\bm x, t)$ with a discrete solution $\hat{\bm u}(\bm x, t)\in \mathbb{R}^m$. As we discussed in the introduction, the discrete solution is represented as the linear sum of $N$ basis functions $\bm \phi_k(\bm x) \in \mathbb{R}^m$ such that the $j$th dimension of $\hat{\bm u}$ is given by $\hat{u}_j = \sum_{k=1}^N c_{jk}(t) \phi_{jk}(\bm x)$. Examples of this discretization process include the finite volume (FV) method, the discontinuous galerkin (DG) method, spectral methods, and the finite element method (FEM). 

In the FV method, the domain $\Omega$ is partitioned into $N$ cells. The solution in each cell is a piecewise constant function; thus the $k$th basis function $\bm \phi_k$ is a vector of ones $\bm 1 \in \mathbb{R}^m$ inside the cell and zero outside the cell. The matrix of coefficients $\bm c_{jk} \in \mathbb{R}^{m \times N}$ thus represents, for each of the $m$ components of $\bm u$, the constant scalar value within each cell. The DG method is similar to the FV method, except that the solution within each cell is represented by a piecewise polynomial of degree $p \in \mathbb{N}$. With DG methods, the domain is  partitioned into cells. Each cell contains a discrete set of basis functions which are usually orthogonal polynomials. A FV method is equivalent to a DG method with $p=0$. Spectral methods represent the solution with global basis functions; often these basis functions are Fourier modes. FEM represents the solution as a sum of polynomial basis functions, though unlike DG methods the solution is piecewise continuous instead of piecewise discontinuous. In this paper we mostly consider FV solvers, though we also consider DG and spectral solvers. We don't consider FEM solvers, nor do we consider other possible basis functions.

Any solver which uses basis functions to represent the solution to \cref{eq:generalpde} will need to find some way of updating the basis function coefficients $\bm c_{jk} \in \mathbb{R}^{m\times N}$ in time. We define a machine learned PDE solver, for the purposes of this paper, as any solver which uses an update rule for the coefficients $\bm c_{jk}$.
This can also be called an `autoregressive' solver. The error-correcting strategy we introduce in \cref{sec:scalarhyperbolic,sec:systemshyperbolic} can be applied to any such solver. Note that the use of ML is not actually part of the above definition.

There are two main types of update rule. One is to use a discrete-time update. The discrete-time update iteratively updates time from $t$ to $t + \Delta t$ and the solution from $\hat{\bm u}$ to $\hat{\bm u} + \Delta t \hat{\mathcal{N}}[\bm \hat{\bm u}]$. The second is to use a continuous-time update. The continuous-time update combined with discretization in space is often called the method of lines (MOL). MOL involves approximating $\hat{\mathcal{N}}$, setting $\frac{\partial \hat{\bm u}}{\partial t} = \hat{\mathcal{N}}$, and using an ODE integrator to advance $\hat{\bm u}$ in time. Using a strong stability preserving Runge Kutta (SSPRK) \cite{ssprk,ssp_gottlieb} ODE integration method and choosing the timestep to satisfy a CFL condition is usually sufficient to ensure that, for hyperbolic PDEs, invariants which are preserved in the continuous-time limit are preserved as time advances. 

Note that some methods which use ML to find a solution to a PDE do not satisfy the above definition of a machine learned PDE solver. For example, the physics-informed neural network (PINN) approach \cite{karniadakis2021physics} does not iteratively update the solution in time, and so the results in this paper do not apply to PINNs. The same is true of the Fourier Neural Operator (FNO) approach \cite{fourier_neural_operator} if the FNO convolves in both space and time. However, the strategy introduced in this paper could be applied to a recurrent FNO which iteratively updates Fourier coefficients in time.

In practical applications, we usually don't care about solving an invariant-preserving equation $\frac{\partial \bm u}{\partial t} + \mathcal{N}[\bm u]=0$, but rather a more complicated equation $\frac{\partial \bm u}{\partial t} + \mathcal{N}[\bm u]=\mathcal{F}[\bm u]$ where $\mathcal{F}[\bm u]$ is some operator that breaks one or more of the invariants preserved by $\mathcal{N}$. In these cases, we can model $\mathcal{N}$ using machine learning, apply an error-correcting algorithm to $\mathcal{N}$, and model $\mathcal{F}$ using some other technique. Doing so ensures that invariants are violated only due to the presence of $\mathcal{F}$, not because of faulty numerics in the calculation of $\mathcal{N}$. Consider a concrete example. The Navier-Stokes equations can be written as a sum of the Euler equations and viscous and forcing terms. The viscous and forcing terms break the invariants of the Euler equations. In this case, we can model the conservative terms in the Euler equations using ML, apply an error-correcting algorithm to these terms, and use some other technique to model the viscous and forcing terms. For an example of how this can be done, see \cite{ml_accelerated_cfd}.

\section{Invariants of Hyperbolic PDEs }\label{sec:standardsolvers}

\noindent PDEs have certain invariants. Computational scientists usually try to design numerical solvers which preserve some, or all, of those invariants. By preserving the right set of invariants, it becomes possible to design solvers which are numerically stable
and which give physically reasonable results when $\Delta x$ and $\Delta t$ are positive. 
In this section, we review the theory of invariant preservation in hyperbolic PDEs. We examine generic scalar hyperbolic PDEs, generic systems of hyperbolic PDEs, and a few example PDEs in each class. We briefly discuss how standard FV methods preserve invariants.


\subsection{Scalar Hyperbolic PDEs}\label{sec:standardsolversscalar}

\noindent Scalar hyperbolic PDEs can be written as
\begin{equation}\label{eq:scalarhyperbolic}
\frac{\partial u}{\partial t} + \bm \nabla \cdot \bm f(u) = 0.
\end{equation}
The solution $u(\bm x, t) \in \mathbb{R}$, $\bm x \in \Omega$, $\Omega \in \mathbb{R}^d$, and the flux $\bm f \in C^1(\mathbb{R}^d)$. 

\noindent \textbf{Continuous Invariants}: Scalar hyperbolic PDEs have one linear invariant which is constant in time and three non-linear invariants which are non-increasing in time \cite{mishra2019numerical}:
\begin{itemize}
    \item Total mass $\int_\Omega u \mathop{d\bm x}$, which is conserved in time.
    \item The $\ell_p$-norm $\int_\Omega |u|^p \mathop{d\bm x}$ for $p > 1$, which is non-increasing in time.
    \item The $\ell_\infty$-norm $\textnormal{max}_{\bm x} u(\bm x, t)$, which is non-increasing in time.
    \item The total variation, which for continuous $u$ in 1D is $\int_0^L \big|\frac{\partial u}{\partial x}\big| \mathop{dx}$. The total variation is non-increasing in time. This is usually called the total variation diminishing (TVD) property.
\end{itemize}
We can prove that hyperbolic PDEs conserve mass by integrating \cref{eq:scalarhyperbolic} over the domain $\Omega$ and using the divergence theorem:
\begin{equation}
    \int_{\Omega} \Big[\frac{\partial u}{\partial t} + \bm \nabla \cdot \bm f(u) \Big]d\bm x = \frac{d}{dt}\int_{\Omega}u \mathop{d\bm x} + \int_{\partial \Omega} \bm f \cdot d\bm n = 0.
\end{equation}
In words: the rate of change of $\int_\Omega u \mathop{d\bm x}$ with respect to time is equal to the negative flux of $u$ through the domain boundary $\partial \Omega$. In an infinite or periodic system $\int_\Omega u \mathop{d\bm x}$ is constant. 

To prove that the $\ell_p$-norm and $\ell_\infty$-norm are non-increasing in time, we use the entropy inequality for scalar hyperbolic PDEs. The entropy inequality can be derived as follows. First multiply \cref{eq:scalarhyperbolic} by a scalar function $\eta'(u)$. $\eta(u)$ is an entropy function. This gives
\begin{equation}\label{eq:entropy_deriv_1}
    \frac{\partial \eta(u)}{\partial t} + \frac{\partial \eta}{\partial u} \Big(\bm \nabla \cdot \bm f(u)\Big) = 0.
\end{equation}
If $u$ is smooth, we can use the chain rule to rearrange \cref{eq:entropy_deriv_1} as
\begin{equation}\label{eq:entropy_deriv_2}
    \frac{\partial \eta}{\partial t} + \frac{\partial \eta}{\partial u} \frac{\partial \bm f}{\partial u} \cdot \bm \nabla u = \frac{\partial \eta}{\partial t} + \frac{\partial \bm \psi}{\partial u} \cdot \bm \nabla u = \frac{\partial \eta}{\partial t} + \bm \nabla \cdot \bm \psi = 0
\end{equation}
where the entropy flux $\bm \psi$ is defined by $\bm \psi'(u) \coloneqq \eta' \bm f'$. \Cref{eq:entropy_deriv_2} states that, for smooth solutions, the entropy $\eta$ satisfies a conservation law. 
For discontinuous $u$, a longer derivation \cite{mishra2019numerical} reveals that the equality in \cref{eq:entropy_deriv_2} is replaced with an inequality for any convex entropy function $\eta(u)$ with corresponding entropy flux $\bm \psi$:
\begin{equation}\label{eq:entropy_inequality}
\frac{\partial \eta(u)}{\partial t} + \bm \nabla \cdot \bm \psi(u) \ge 0.
\end{equation}
Integrating \cref{eq:entropy_inequality} over $\Omega$ shows that the rate of change of total entropy $\int_\Omega \eta \mathop{d\bm x}$ is equal to the negative entropy flux through the domain boundary $\partial \Omega$; in an infinite or periodic system total entropy is non-decreasing.
In brief: for continuous $u$ entropy is conserved, while for discontinuous $u$ entropy increases.

By choosing $\eta(u) = -|u|^p$ and integrating \cref{eq:entropy_inequality} over $\Omega$, we have the non-increasing $\ell_p$-norm property. Taking the limit as $p \rightarrow \infty$ gives the non-increasing $\ell_\infty$-norm property. The TVD property is derived in \cite{mishra2019numerical}.

\textbf{Finite volume (FV) method}: A common approach for solving hyperbolic PDEs is by using a finite volume (FV) method. FV methods divide the spatial domain  $\Omega$ into a number of discrete cells $\Omega_j$, then use a scalar value to represent the solution average within each cell. For example, on the 1D domain $x \in [0,L]$ with uniform cell width, a FV method divides the domain into $N$ cells of width $\Delta x = \nicefrac{L}{N}$ where the left and right boundaries of the $j$th cell for $j = 1, \dots, N$ are $x_{j-\nicefrac{1}{2}} = (j-1)\Delta x$ and $x_{j+\nicefrac{1}{2}}=j\Delta x$ respectively. FV methods use a scalar value $u_j(t)$ to represent the solution average within each cell where $u_j(t) \coloneqq \int_{x_{j-\nicefrac{1}{2}}}^{x_{j+\nicefrac{1}{2}}} u(x, t)  \mathop{dx}$. The standard FV equations for the time-derivative of $u_j$ in 1D and $u_{i,j}$ in 2D are simply discrete versions of \cref{eq:scalarhyperbolic}:
\begin{subequations}
\begin{equation}\label{eq:fv_a}
    \frac{\partial u_j}{\partial t} + \frac{f_{j+\frac{1}{2}} - f_{j-\frac{1}{2}}}{\Delta x} = 0
\end{equation}
\begin{equation}\label{eq:fv_b}
    \frac{\partial u_{i,j}}{\partial t} + \frac{f^x_{i+\frac{1}{2},j} - f^x_{i-\frac{1}{2},j}}{\Delta x} + \frac{f^y_{i,j+\frac{1}{2}} - f^y_{i,j-\frac{1}{2}}}{\Delta y} = 0.
\end{equation}
$f_{j+\nicefrac{1}{2}}$ is the flux at the cell boundary $x_{j+\nicefrac{1}{2}}$. $f^x_{i+\nicefrac{1}{2},j}$ and $f^y_{i,j+\nicefrac{1}{2}}$ are the average x-directed and y-directed fluxes through the right and top cell boundaries, e.g., $f^x_{i+\nicefrac{1}{2},j} \coloneqq \frac{1}{\Delta y}\int_{y=y_{j-\nicefrac{1}{2}}}^{y=y_{j+\nicefrac{1}{2}}} \bm{\hat{x}} \cdot \bm f(x_{i+\nicefrac{1}{2}},y) \mathop{dy}$. 
In higher dimensions, the FV update equation in the $j$th cell can be written as
\begin{equation}\label{eq:fv_c}
    \frac{\partial u_j}{\partial t} + \frac{1}{|\Omega_j|}\oint_{\partial \Omega_j} \bm{f}\cdot \mathop{d\bm s} = 0
\end{equation}
\end{subequations}
where $|\Omega_j| \coloneqq \int_{\Omega_j} \mathop{d\bm x}$ is the volume in cell $\Omega_j$ and $\mathop{d\bm s}$ is the outward normal vector at cell boundary $\partial \Omega_j$.
In 1D, \cref{eq:fv_a} can be derived by applying the integral $\int_{x_{j-\nicefrac{1}{2}}}^{x_{j+\nicefrac{1}{2}}} (...) \mathop{dx}$ to \cref{eq:scalarhyperbolic} for all $j \in 1, \dots, N$; a similar calculation gives \cref{eq:fv_b} and \cref{eq:fv_c}.
So long as $f_{j+\nicefrac{1}{2}}$ or $f_{i+\nicefrac{1}{2},j}^x$ and $f_{i,j+\nicefrac{1}{2}}^y$ are exact for all $t$, then $u_j$ or $u_{ij}$ will be exact for all $t$. Thus, the key challenge for a FV scheme is to accurately reconstruct the flux at cell boundaries.

\noindent \textbf{Discrete Invariants}: FV schemes conserve a discrete analogue of the continuous linear invariant $\int_\Omega u \mathop{d \bm x}$ by construction. In 1D, we can see this with a short proof:  
$\nicefrac{d}{d t} \sum_{j=1}^N u_j \Delta x = \Delta x \sum_{j=1}^N \nicefrac{\partial u_j}{\partial t} = - \sum_{j=1}^N (f_{j+\nicefrac{1}{2}} - f_{j-\nicefrac{1}{2}}) = f_{\nicefrac{1}{2}} - f_{N+\nicefrac{1}{2}}$. The rate of change of the discrete mass is equal to the flux of $u$ through the boundaries; in a periodic system this equals 0.

Although FV schemes preserve a discrete analogue of conservation of mass by construction, they do not automatically preserve discrete analogues of any of the non-linear invariants of the continuous PDE. Instead, FV methods preserve non-linear invariants through careful choice of flux.

The only known way of inheriting discrete analogues of all three non-linear invariants of \cref{eq:scalarhyperbolic} (non-increasing $\ell_p$-norm, non-increasing $\ell_\infty$-norm, and TVD) is to use a consistent monotone flux function while satisfying a CFL condition \cite{mishra2019numerical}. See \cite{dale_durran} for definitions of consistency and monotonicity. An example of a monotone flux function for the linear advection equation $f=cu$ is the upwind flux
\begin{equation}
    f_{j+\frac{1}{2}} = \begin{cases}
        c u_j &\textnormal{ if } c \ge 0 \\
        c u_{j+1} &\textnormal{ if } c < 0.
    \end{cases}
\end{equation}
For non-linear $f(u)$, a Reimann solver or approximate Reimann solver results in a monotone flux function.  Examples of monotone flux functions include the Godunov flux 
\begin{equation}
    f_{j+\frac{1}{2}} = \begin{cases}
        \textnormal{min}_{u_j \le u \le u_{j+1}} f(u) &\textnormal{ if } u_j \le u_{j+1} \\
        \textnormal{max}_{u_j \le u \le u_{j+1}} f(u) &\textnormal{ if } u_j > u_{j+1}
    \end{cases}
\end{equation}
and the Lax-Friedrichs flux
\begin{equation}
    f_{j+\frac{1}{2}} = \frac{f(u_j) + f(u_{j+1})}{2} - \frac{\Delta x}{2\Delta t} (u_{j+1} - u_j).
\end{equation}
Unfortunately, Godunov's famous theorem from 1959 implies that monotone schemes can be at most first-order accurate \cite{godunov1959finite}. This means that while monotone schemes preserve all the invariants of the underlying PDE, they are usually not very accurate.

Godunov's theorem reveals a more general lesson. For some PDEs, it is impossible to design highly accurate numerical methods that preserve discrete analogues of \textit{every} invariant of the continuous PDE. Designers of numerical methods must therefore determine which invariants of the continuous system should be preserved by the discrete system and which invariants either cannot be preserved or degrade the accuracy of the discrete system.

For scalar hyperbolic PDEs, it turns out that it is possible to design accurate and stable numerical solvers by preserving just \textit{one} of the three non-linear invariants of \cref{eq:scalarhyperbolic} \cite{dale_durran}.
One such scheme is the MUSCL scheme, introduced in a seminal paper by Van Leer \cite{muscl}. MUSCL uses limiters to reconstruct the solution at cell boundaries. In 1D, this preserves a discrete analogue of the TVD property. Doing so guarantees numerical stability and prevents spurious oscillations\footnote{Spurious oscillations are unphysical oscillations which develop in numerical methods that do not have enough numerical diffusion to damp high-$k$ modes that develop near steep gradients \cite{john2007spurious,shocks_artificial_viscosity}.},  while retaining second-order accuracy.

It is also possible to design stable numerical methods by preserving a discrete analogue of an $\ell_p$-norm non-increasing property of the solution. For the linear advection equation $f=cu$, the centered flux 
\begin{equation}
    f_{j+\frac{1}{2}} = \frac{c}{2}(u_j + u_{j+1})
\end{equation}
conserves the discrete $\ell_2$-norm in the continuous-time limit, though if used with a forward Euler update the centered flux increases the discrete $\ell_2$-norm leading to numerical instability. For non-linear $f(u)$ the flux formula
\begin{equation}
    f_{j+\frac{1}{2}} = \int_0^1 f(\hat{u})d\theta \textnormal{\hspace{0.1cm} where \hspace{0.1cm}} \hat{u}(\theta) = u_j + \theta (u_{j+1}-u_j)
\end{equation}
conserves the discrete $\ell_2$-norm \cite{jameson2008construction} in the continuous-time limit.

\subsubsection{2D Incompressible Euler Equations\label{sec:standardsolversincompressible}}

\noindent The 2D incompressible Euler equations in vorticity form are a scalar hyperbolic PDE coupled to an elliptic PDE:
\begin{align}\label{eq:euler}
    \frac{\partial \chi}{\partial t} + \bm \nabla \cdot (\bm u \chi) = 0\textnormal{,} && \bm u = \bm \nabla \psi \times \hat{e}_z\textnormal{,} && -\bm\nabla^2 \psi = \chi.
\end{align}
These equations can be written as a Hamiltonian system (\cref{eq:hamiltoniansystem}) for the vorticity $\chi(x,y,t)$ evolving under a Hamiltonian given by the streamfunction $\psi(x,y,t)$ \cite{hamiltonian_structure}. 

\noindent \textbf{Continuous Invariants}: \Cref{eq:euler} has the same invariants as the generic scalar hyperbolic \cref{eq:scalarhyperbolic}. In particular, \cref{eq:euler} conserves the mass $\int_\Omega \chi \mathop{dx}\mathop{dy}$ and the $\ell_2$-norm $\int_{\Omega} \chi^2 \mathop{dx}\mathop{dy}$, also called the Enstrophy. \Cref{eq:euler} has an additional conserved invariant, the energy $\frac{1}{2} \int_{\Omega} \bm u^2 \mathop{dx}\mathop{dy}$.
We introduce the notation $\mathop{dz} = \mathop{dx}\mathop{dy}$. We prove conservation of mass using integration by parts, and canceling the boundary term using periodic BCs:
\begin{equation}
    \frac{d}{dt} \int_\Omega \chi \mathop{d z}= \int_\Omega \frac{\partial \chi}{\partial t} \mathop{d z} = - \int_\Omega \bm \nabla \cdot (\bm u \chi) \mathop{d z} = - \int_{\partial \Omega} \bm u\chi \cdot \mathop{d\bm s} = 0.
\end{equation}
We prove Enstrophy conservation as follows. From incompressibility, $\bm \nabla \cdot \bm u = 0$, which implies $\bm \nabla \cdot (\bm u \nicefrac{\chi^2}{2}) =\chi\bm u\cdot \bm \nabla \chi = \chi \bm \nabla \cdot (\bm u \chi)$. Using Gauss's theorem and periodicity,
\begin{equation}
    \frac{d}{dt} \int_\Omega \frac{1}{2}\chi^2 \mathop{d z} = - \int_\Omega \chi \bm \nabla \cdot (\bm u \chi) \mathop{d z} = -\int_{\partial \Omega} \frac{1}{2} {\bm u \chi^2} \cdot \mathop{d\bm s} = 0.
\end{equation}
We prove energy conservation using integration by parts and $\bm u \cdot \bm \nabla \psi = 0$:
\begin{equation}
\begin{split}
    \frac{d}{dt}\int_{\Omega} \frac{1}{2}\bm u^2 \mathop{d z} &= \frac{d}{dt}\int_{\Omega} \frac{1}{2}(\bm \nabla \psi)^2 \mathop{d z} = \int_{\Omega} \bm \nabla \psi \cdot \bm \nabla \frac{\partial \psi}{\partial t} \mathop{d z}\\ &=- \int_\Omega \psi \nabla^2 \frac{\partial \psi}{\partial t} \mathop{d z} = \int_\Omega \psi \frac{\partial \chi}{\partial t} \mathop{d z} \\&= -\int_\Omega \psi \bm \nabla \cdot (\bm u \chi) = \int_\Omega \chi \bm u \cdot \bm \nabla \psi \mathop{d z} = 0.
\end{split}
\end{equation}

\subsection{Systems of Hyberbolic PDEs\label{sec:standard_systems}}

\noindent We consider systems of hyperbolic PDEs in 1D written in conservation form, given by
\begin{equation}\label{eq:system_conservation_form}
    \frac{\partial \bm u}{\partial t} + \frac{\partial}{\partial x} \bm F(\bm u) = 0
\end{equation}
where $\bm u \in \mathbb{R}^m$ and $\bm F(u) \in C^1(\mathbb{R}^m)$.
\Cref{eq:system_conservation_form} is hyperbolic if the Jacobian matrix $\frac{\partial \bm F}{\partial \bm u}$ has real eigenvalues and a complete set of linearly independent eigenvectors \cite{leveque_green_book}.

\noindent \textbf{Continuous Invariants}: \Cref{eq:system_conservation_form} implies that each component of $\int \bm u(x, t) \mathop{dx}$ is conserved. In a 1D periodic system with $x \in [0, L]$, an integral over $x$ makes this apparent: $\frac{d}{dt}\int_{0}^L \bm u \mathop{dx} = \int \frac{\partial \bm u}{\partial t} \mathop{dx} = - \int_0^L \frac{\partial \bm F}{\partial x} \mathop{dx} = \bm F(0) - \bm F(L)$. 
The total rate of change of $\bm u$ equals the flux $\bm F$ through the boundaries; in a periodic or infinite system this equals zero.

For certain systems of hyperbolic PDEs, it is possible to define a generalized scalar convex entropy function $\eta(\bm u)$ and entropy flux $\psi(\bm u)$ such that the entropy satisfies an inequality \cite{leveque_green_book}
\begin{equation}\label{eq:system_entropy_inequality}
    \frac{\partial \eta(\bm u)}{\partial t} + \frac{\partial\psi(\bm u)}{\partial x}  \ge 0.
\end{equation} 
Integrating \cref{eq:system_entropy_inequality} over $x$ for a 1D periodic system where $x \in [0, L]$ shows that the total entropy is non-decreasing: $\frac{d}{dt} \int \eta(\bm u) \mathop{dx} \ge \psi(0) - \psi(L) = 0$.
It can be shown that \cref{eq:system_conservation_form} satisfies the entropy inequality \cref{eq:system_entropy_inequality} if there exists a change of variables from $\bm u$ to $\bm w$ such that $\nicefrac{\partial \bm u}{\partial \bm w}$ and $\nicefrac{\partial \bm F}{\partial \bm w}$ are both symmetric \cite{jameson2008construction,harten1983symmetric}. 
If this change of variables is made, it turns out that the entropy variable $\bm w = \nicefrac{\partial \eta}{\partial \bm u}$, $\nicefrac{\partial \eta}{\partial t} = \bm w^T \nicefrac{\partial \bm u}{\partial t}$, and $\nicefrac{\partial \psi}{\partial \bm u} = (\nicefrac{\partial \eta}{\partial \bm u})^T \nicefrac{\partial \bm F}{\partial \bm u}$.

Notice that while generic scalar hyperbolic PDEs all conserve the same non-linear invariants, generic systems of hyperbolic PDEs are not guaranteed to have any non-linear invariants. Many physically relevant systems of hyperbolic PDEs do have non-linear invariants, but these invariants differ between systems of PDEs.

\textbf{Discrete Invariants}: FV schemes conserve a discrete analogue of the continuous linear invariant $\int_\Omega \bm u \mathop{d x}$ by construction. In 1D, we can see this with a short proof:  
$\nicefrac{d}{d t} \sum_{j=1}^N \bm u_j \Delta x = \Delta x \sum_{j=1}^N \nicefrac{\partial \bm u_j}{\partial t} = - \sum_{j=1}^N (\bm F_{j+\nicefrac{1}{2}} - \bm F_{j-\nicefrac{1}{2}}) = \bm F_{N+\nicefrac{1}{2}}- \bm F_{\nicefrac{1}{2}}$. The rate of change of the discrete mass is equal to the flux of $\bm u$ through the boundaries; in a periodic system this equals 0.

When a generalized entropy function $\eta(\bm u)$ exists, FV schemes do not guarantee that the discrete entropy $\sum_{j=1}^N \eta_j(\bm u_j) \Delta x$ is non-decreasing in time. Instead, FV schemes inherit a discrete analogue of the non-decreasing total entropy property by satisfying a discrete entropy inequality in each grid cell:
\begin{equation}\label{eq:discrete_entropy_inequality}
    \frac{\partial \eta_j}{\partial t} + \frac{\psi_{j+\frac{1}{2}} - \psi_{j - \frac{1}{2}}}{\Delta x} \ge 0.
\end{equation}
One way of satisfying \cref{eq:discrete_entropy_inequality} is to use a Reimann solver or approximate Reimann solver to compute $\bm F_{j+\nicefrac{1}{2}}$ \cite{leveque_green_book}. Doing so preserves the non-decreasing total entropy invariant.

\subsubsection{Compressible Euler Equations}

\noindent The compressible Euler equations of gas dynamics in 1D are given by 
\begin{equation}\label{eq:compressibleeuler}
    \frac{\partial }{\partial t} \begin{bmatrix}
           \rho \\
           \rho v \\
           E
         \end{bmatrix} + \frac{\partial}{\partial x} \begin{bmatrix}
             \rho v \\
             \rho v^2 + p \\
            v(E+p)
         \end{bmatrix} = 0.
\end{equation}
where $\rho$ is the density, $u$ is the velocity, $p$ is the pressure and $E$ is the energy. The equation of state for an ideal gas is
\begin{equation}\label{eq:equationofstate}
    E = \frac{p}{\gamma - 1} + \frac{1}{2}\rho v^2
\end{equation}
where $\gamma$ is the ratio of specific heat at constant pressure to the specific heat at constant volume. In the notation of \cref{eq:system_conservation_form}, $\bm u = \begin{bmatrix}
    \rho, \rho v, E
\end{bmatrix}$ and $\bm F = \begin{bmatrix}
    \rho v, \rho v^2 + p, v(E+p)
\end{bmatrix}$.

\noindent \textbf{Continuous Invariants}: As with all hyperbolic PDEs, the time rate of change of the linear invariant $\int \bm u \mathop{dx}$ is equal to the negative flux $\bm F$ of $\bm u$ through the domain boundaries. In an infinite or periodic system $\int \bm u \mathop{dx}$ is constant in time. 

The Euler equations have two positivity invariants: the density $\rho \in \mathbb{R}$ and the pressure $p \in \mathbb{R}$ are both everywhere non-negative. From the equation of state \cref{eq:equationofstate}, non-negativity of density and pressure imply that the energy $E \in \mathbb{R}$ is everywhere non-negative.

The Euler equations also satisfy an entropy inequality \cref{eq:system_entropy_inequality}. This implies that the total entropy $\int \eta(\bm u) \mathop{dx}$ is non-decreasing in time. It has been shown \cite{harten1983symmetric} that the generalized entropy function $\eta(s) = \rho g(s)$ is convex for any $g(s)$ with specific entropy $s = \log{(\nicefrac{p}{\rho^\gamma})}$ for which $\nicefrac{g''}{g'} < \gamma^{-1}$. In this case the entropy flux $\psi = \rho g(s) v$. 
If we choose $g(s) = e^{\nicefrac{s}{\gamma + 1}}$, some tedious algebra shows that the entropy variable $\bm w = \nicefrac{\partial \eta}{\partial \bm u}$ equals
\begin{equation}\label{eq:compressibleeuler_entropy_variable}
    \bm w = \frac{p^*}{p}\begin{bmatrix}
        E \\
        -\rho v \\
        \rho
    \end{bmatrix}
\end{equation}
where
\begin{equation}
    p^* = \frac{\gamma - 1}{\gamma + 1} \bigg(\frac{p}{\rho^\gamma}\bigg)^{\frac{1}{\gamma + 1}}.
\end{equation}

\section{Invariant-Preserving Error-Correcting Algorithms for Scalar Hyperbolic PDEs}\label{sec:scalarhyperbolic}

\noindent We now turn to the main purpose of this paper: designing machine learned solvers that preserve discrete invariants.
As we discussed in \cref{sec:standardsolversscalar}, scalar hyperbolic PDEs conserve the linear invariant $\int_\Omega u \mathop{d\bm x}$, have non-increasing $\ell_p$-norm, have non-increasing $\ell_\infty$-norm, and are TVD. The only known way of inheriting discrete analogues of all three non-linear invariants of \cref{eq:scalarhyperbolic} is to use a consistent monotone flux function. Godunov's theorem implies that such schemes can be at most first-order accurate. Thus, we will not try to design monotone machine learned solvers. Instead, we will design machine learned solvers that preserve two invariants: mass conservation and one of the non-linear invariants. This ensures that the machine learned solver is numerically stable, while being flexible enough to accurately predict the rate of change of the solution. 

In this section, we introduce error-correcting algorithms for scalar hyperbolic PDEs in periodic domains which modify an update rule to enforce mass conservation and change the time-derivative of the discrete $\ell_2$-norm from $\nicefrac{d\ell_2^{\textnormal{old}}}{dt}$ to $\nicefrac{d\ell_2^{\textnormal{new}}}{dt}$.
These error-correcting algorithms make it easy to design invariant-preserving machine learned solvers: first, at each timestep use the error-correcting algorithm to modify the output of the update rule to ensure that the discrete mass is conserved. Second, if $\nicefrac{d\ell_2^{\textnormal{old}}}{dt} > 0$, set $\nicefrac{d\ell_2^{\textnormal{new}}}{dt} =0$.
If $\nicefrac{d\ell_2^{\textnormal{old}}}{dt} \le 0$, set $\nicefrac{d\ell_2^{\textnormal{new}}}{dt} = \nicefrac{d\ell_2^{\textnormal{old}}}{dt}$. If we have a priori information about the expected rate of change of the discrete $\ell_2$-norm, we can set $\nicefrac{d\ell_2^{\textnormal{new}}}{dt}$ to that value.

We consider a variety of different solver types and update rules. \Cref{sec:continuoustimeFV,sec:fluxpredictingFV,sec:DG,sec:fourier} show how to design invariant-preserving FV solvers, DG solvers, and spectral solvers with a continuous-time update. \Cref{sec:continuoustimeFV} is derived for arbitrary mesh shapes. \Cref{sec:discretetimeFV} shows how to design invariant-preserving FV solvers when using a discrete-time update.
In \cref{sec:energyconservation} we consider a different scalar hyperbolic PDE, the 2D incompressible Euler equations. For these equations we show how to preserve an additional invariant, conservation of energy.
For each type of solver, the general strategy is the same: at each timestep, apply an error-correcting algorithm to the update rule.
We only consider periodic BCs in this section. With the exception of the spectral solver in \cref{sec:fourier}, we can derive invariant-preserving algorithms for non-periodic BCs by estimating the flux through the boundaries.

\subsection{Flux-Predicting FV Solvers \label{sec:fluxpredictingFV}}

\noindent Suppose we are interested in designing a machine learned solver for \cref{eq:scalarhyperbolic} which, like the FV method, divides the domain $\Omega$ into a number of grid cells $\Omega_j$ and represents the solution average in the $j$th cell as a scalar $u_j$. Suppose also we use the continuous-time FV update \cref{eq:fv_a,eq:fv_b,eq:fv_c} to compute $\nicefrac{\partial u_j}{\partial t}$. Because \cref{eq:fv_a,eq:fv_b,eq:fv_c} are exact, to evolve the discrete solution accurately a machine learned solver has to predict the (average) flux across cell boundaries accurately. We refer to this type of machine learned solver as a `flux-predicting finite volume' solver. 

Flux-predicting FV solvers conserve the discrete mass by construction. We now show how to modify the predicted flux to ensure that the discrete $\ell_2$-norm is non-increasing in the continuous-time limit. In 1D, the discrete $\ell_2$-norm is non-increasing if 
\begin{equation}
    \frac{d}{dt}\sum_{j=1}^N \frac{1}{2}u_j^2 \Delta x_j \le 0
\end{equation}
for all $t$. $\Delta x_j$ is the width of cell $\Omega_j$. Some simple algebra and \cref{eq:fv_a} gives
\begin{equation}
    \frac{d}{dt} \sum_{j=1}^N \frac{\Delta x_j}{2} u_j^2 = \sum_{j=1}^N u_j \frac{\partial u_j} {\partial t} \Delta x_j =  -\sum_{j=1}^N u_j (f_{j+\frac{1}{2}} - f_{j-\frac{1}{2}}).
\end{equation}
Performing summation by parts gives
\begin{equation}
-\sum_{j=1}^N u_j (f_{j+\frac{1}{2}} - f_{j-\frac{1}{2}}) = \sum_{j=1}^{N-1} f_{j+\frac{1}{2}}\big(u_{j+1} - u_j \big) + f_{\frac{1}{2}} u_1 - f_{N+\frac{1}{2}} u_N \le 0.
\end{equation}
In a periodic domain, this is simply
\begin{equation}\label{eq:stability_energy_method}
    \frac{d}{dt}\sum_{j=1}^N  \frac{\Delta x_j }{2}(u_j)^2 = \sum_{j=1}^{N} f_{j+\frac{1}{2}}\big(u_{j+1} - u_j \big) \le 0.
\end{equation}
For the rest of the section we will assume a periodic domain, though our results can easily be generalized to non-periodic domains. 
Let us now define 
\begin{equation}
    \frac{d\ell_{2}^{\textnormal{old}}}{dt} \coloneqq \sum_{j=1}^N f_{j+\frac{1}{2}}(u_{j+1} - u_j) 
\end{equation}
as the original rate of change of the discrete $\ell_2$-norm, and $\nicefrac{d\ell_{2}^{\textnormal{new}}}{dt}$ as the desired rate of change of the discrete $\ell_2$-norm. To ensure non-increasing $\ell_2$-norm, we want $\nicefrac{d\ell_2^{\textnormal{new}}}{dt}\le 0$.
We also define $\bm u_j \coloneqq \{u_j\}_{j=1}^N$ as a vector representation of the discrete solution. 
We can change the time-derivative of the discrete $\ell_2$-norm from $\nicefrac{d\ell_2^{\textnormal{old}}}{dt}$ to $\nicefrac{d\ell_2^{\textnormal{new}}}{dt}$ by making the following transformation to $f_{j+\nicefrac{1}{2}}$:
\begin{equation}\label{eq:1d_stability}
	f_{j+\frac{1}{2}} \Rightarrow 
	f_{j+\frac{1}{2}} + 
	\frac{(\nicefrac{d\ell_2^{\textnormal{new}}}{dt} - \nicefrac{d\ell_2^{\textnormal{old}}}{dt}) G_{j+\nicefrac{1}{2}}(\bm u_j)}
		{\sum_{k=1}^N G_{k+\nicefrac{1}{2}}(\bm u_k)(u_{k+1} - u_k) }
\end{equation}
for any scalar $\nicefrac{d\ell_2^{\textnormal{new}}}{dt}$ and any non-constant, finite function $G_{j+\nicefrac{1}{2}}(\bm u_j)$ for which 
\begin{equation}
    \sum_{k=1}^N G_{k+\nicefrac{1}{2}}(\bm u_k)(u_{k+1} - u_k) \ne 0.
\end{equation}
As the reader can verify by plugging \cref{eq:1d_stability} into \cref{eq:stability_energy_method}, \cref{eq:1d_stability} modifies $f_{j+\nicefrac{1}{2}}$ in a way that adds a constant $(\nicefrac{d\ell_2^{\textnormal{new}}}{dt} - \nicefrac{d\ell_2^{\textnormal{old}}}{dt})$ to \cref{eq:stability_energy_method} via cancellation of the denominator.
Note that $G_{j+\nicefrac{1}{2}}(\bm u_j)$ is a hyperparameter that determines how each $f_{j+\nicefrac{1}{2}}$ is modified and $\nicefrac{d\ell_2^{\textnormal{new}}}{dt}$ is a user-defined quantity which sets the rate of change of the discrete $\ell_2$-norm. 
A similar calculation in a 2D periodic rectangular domain with uniform grid spacing reveals that the rate of change of the discrete $\ell_2$-norm is given by
\begin{equation}\label{eq:l2cons2d}
\begin{split}
    \frac{d}{dt}\sum_{i,j} \frac{u_{i,j}^2}{2} \Delta x \Delta y = &\Delta y\sum_{i,j} f^x_{i+\frac{1}{2},j}(u_{i+1,j} - u_{i,j}) \\ + &\Delta x \sum_{i,j} f^y_{i,j+\frac{1}{2}}(u_{i,j+1} - u_{i,j}) \le 0.
\end{split}
\end{equation}
We define 
\begin{subequations}
\begin{equation}
    \frac{d\ell^{\textnormal{old},x}_{2}}{dt} \coloneqq \Delta y \sum_{i,j} f^x_{i+\frac{1}{2},j}(u_{i+1,j} - u_{i,j}),
\end{equation}
\begin{equation}
    \frac{d\ell^{\textnormal{old},y}_{2}}{dt} \coloneqq \Delta x \sum_{i,j} f^y_{i,j+\frac{1}{2}} (u_{i,j+1} - u_{i,j}.
\end{equation}
\end{subequations}
\Cref{eq:l2cons2d} will be satisfied if the following transformations are made to $f^x_{i+\frac{1}{2},j}$ and $f^y_{i,j+\frac{1}{2}}$:
\begin{subequations}
\begin{equation} \label{eq:2d_stability_a}
	f^x_{i+\frac{1}{2},j} \Rightarrow 
	f^x_{i+\frac{1}{2},j} + \frac{(\nicefrac{d\ell^{\textnormal{new},x}_{2}}{dt} - \nicefrac{d\ell^{\textnormal{old},x}_{2}}{dt}) G^x_{i+\nicefrac{1}{2},j}(\bm u_{ij})}{\Delta y\sum_{k,l} G^x_{k+\nicefrac{1}{2},l}(\bm u_{kl})(u_{k+1,l} - u_{k,l})}
\end{equation}
\begin{equation} \label{eq:2d_stability_b} 
	f^y_{i,j+\frac{1}{2}} \Rightarrow 
	f^y_{i,j+\frac{1}{2}} + \frac{(\nicefrac{d\ell^{\textnormal{new},y}_{2}}{dt} - \nicefrac{d\ell^{\textnormal{old},y}_{2}}{dt}) G^y_{i,j+\nicefrac{1}{2}}(\bm u_{ij})}{\Delta x\sum_{k,l} G^y_{k,l+\nicefrac{1}{2}}(\bm u_{kl}) (u_{k,l+1} - u_{k,l})}
\end{equation}
\end{subequations}
for any scalars $\nicefrac{d\ell^{\textnormal{new},x}_{2}}{dt}$ and $\nicefrac{d\ell^{\textnormal{new},y}_{2}}{dt}$ where $\nicefrac{d\ell^{\textnormal{new},x}_{2}}{dt} + \nicefrac{d\ell^{\textnormal{new},y}_{2}}{dt} \le 0$ and any non-constant, finite functions $G^x_{i+\nicefrac{1}{2},j}(\bm u_{ij})$ and $G^y_{i,j+\nicefrac{1}{2}}(\bm u_{ij})$ for which 
\begin{subequations}
    \begin{equation}
        \sum_{k,l} G^x_{k+\nicefrac{1}{2},l}(\bm u_{kl})(u_{k+1,l} - u_{k,l}) \ne 0
    \end{equation}
    \begin{equation}
        \sum_{k,l} G^y_{k,l+\nicefrac{1}{2}}(\bm u_{kl}) (u_{k,l+1} - u_{k,l}) \ne 0
    \end{equation}
\end{subequations}
\Cref{eq:1d_stability,eq:2d_stability_a,eq:2d_stability_b} are the main results of section \ref{sec:fluxpredictingFV}; for scalar hyperbolic PDEs in 1D and 2D periodic domains they ensure that the discrete $\ell_2$-norm will be non-increasing in the continuous-time limit.

How should the hyperparameters $G_{j+\nicefrac{1}{2}}(\bm u_j)$, $G^x_{i+\nicefrac{1}{2},j}(\bm u_{ij})$, and $G^y_{i,j+\nicefrac{1}{2}}(\bm u_{ij})$ be set? In our experiments, we set $G_{j+\nicefrac{1}{2}}(\bm u_j) = (u_{j+1}-u_j)$, $G^x_{i+\nicefrac{1}{2},j}(\bm u_{ij}) = (u_{i+1,j} - u_{i,j})$, and $G^y_{i,j+\nicefrac{1}{2}}(\bm u_{ij}) = (u_{i,j+1} - u_{i,j})$. These choices have a simple physical interpretation: they correspond to the addition of a spatially constant diffusion coefficient everywhere in space \cite{artificial_viscosity}. Possible alternatives include setting $G_{j+\nicefrac{1}{2}}(\bm u_j) = (u_{j+1}-u_j)^\beta$ for $\beta>1$ or $G_{j+\nicefrac{1}{2}}(\bm u_j) = \alpha_{j+\nicefrac{1}{2}} (u_{j+1}-u_j)$ for $\alpha_{j+\nicefrac{1}{2}} \in \mathbb{R}$. Choosing large $\beta$ increases the amount of numerical diffusion added at discontinuities and decreases the amount of diffusion added in smooth regions, while $\alpha_{j+\nicefrac{1}{2}}$ is a spatially dependent scalar which determines a spatially varying distribution of added numerical diffusion.

\begin{figure}
  \centering
  \includegraphics[width=\textwidth]{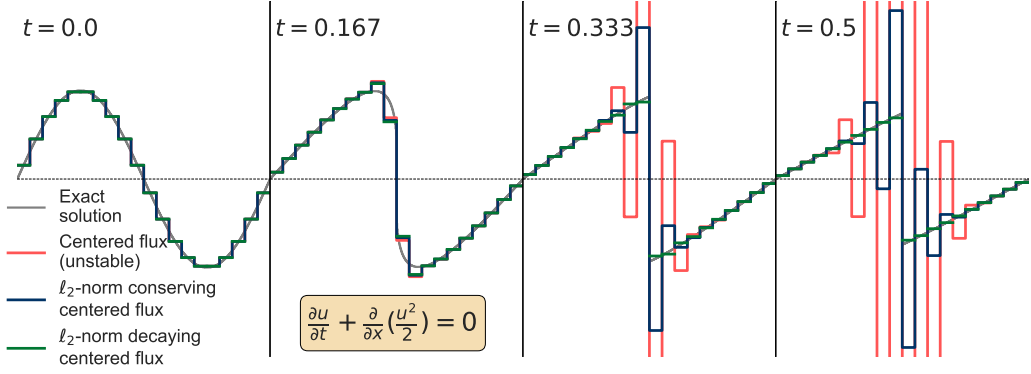}
  \caption{We modify the flux at cell boundaries to ensure non-increasing $\ell_2$-norm, thereby turning an unstable solver into a stable solver. While the centered flux $f_{j+\nicefrac{1}{2}}=\nicefrac{(u_{j}^2 + u_{j+1}^2)}{4}$ (red) is an unstable choice of flux on the inviscid Burgers equation and blows up by $t=0.5$, \cref{eq:1d_stability} with $\nicefrac{d\ell_2^{\textnormal{new}}}{dt}=0$ (blue) ensures that the discrete $\ell_2$-norm is conserved and is thus stable. Setting $\nicefrac{d\ell_2^{\textnormal{new}}}{dt} = \nicefrac{d\ell_2^{\textnormal{exact}}}{dt}$ (green) results in a more accurate solver.}
  \label{fig:fluxFV}
\end{figure}

We now illustrate the effect of modifying an unstable solver in 1D using \cref{eq:1d_stability}.
In \cref{fig:fluxFV}, we solve the inviscid Burgers' equation using a 3rd-order SSPRK ODE integrator \cite{ssprk} with the initial condition $u_0(x) = \sin{x}$.
It turns out that the centered flux $f_{j+\nicefrac{1}{2}}=\nicefrac{(u_j^2 + u_{j+1}^2)}{4}$, shown in red in \cref{fig:fluxFV}, gives $\nicefrac{d\ell_2^{\textnormal{old}}}{dt} > 0$ and is thus unstable and inaccurate.
If we transform $f_{j+\nicefrac{1}{2}}$ according to \cref{eq:1d_stability} with $\nicefrac{d\ell_2^{\textnormal{new}}}{dt}=0$, shown in blue, the solver becomes $\ell_2$-norm conserving and is stable.
Because the $\ell_2$-norm of the exact solution (grey) is decreasing in time, we can further improve the accuracy of the unstable centered flux by setting $\nicefrac{d\ell_2^{\textnormal{new}}}{dt} < 0$.
If we instead transform $f_{j+\nicefrac{1}{2}}$ according to \cref{eq:1d_stability} with $\nicefrac{d\ell_2^{\textnormal{new}}}{dt}=\nicefrac{d\ell_2^{\textnormal{exact}}}{dt}$, shown in green, the solver becomes much more accurate.

Note that the blue solution in \cref{fig:fluxFV} conserves both the discrete mass and the discrete $\ell_2$-norm, but does not maintain a discrete analogue of the total variation diminishing (TVD) property \cite{dale_durran} of the scalar Burgers equation. As a result, the blue solution does not fully eliminate the high-$k$ oscillations that develop in the red solution.
Note also that, in practice, we will usually not have a priori knowledge of $\nicefrac{d\ell_2^{\textnormal{exact}}}{dt}$.
Even if we did, it would be unrealistic to expect that error-correcting invariant preservation algorithms will in general be able to turn a very inaccurate solver (red in \cref{fig:fluxFV}) into an accurate solver (green in \cref{fig:fluxFV}). The takeaway from \cref{fig:fluxFV} is simply that solvers which preserve the right set of invariants tend to be more accurate than solvers which do not preserve those invariants.

\subsection{Continuous-Time FV Solvers with Arbitrary Time-Derivative \label{sec:continuoustimeFV}}

\noindent In section \ref{sec:fluxpredictingFV}, we considered schemes that use ML to predict the flux $f$ across cell boundaries. Using integration by parts, we found that we could modify the fluxes to control the rate of change of the discrete $\ell_2$-norm and thereby preserve the non-increasing $\ell_2$-norm invariant. However, some machine learned PDE solvers may use an alternative form for the time-derivative which does not involve predicting the flux across cell boundaries. Thus, we now consider the more general problem of how to design invariant-preserving solvers for \cref{eq:scalarhyperbolic} with arbitrary time-derivative in arbitrary number of dimensions with arbitrary cell shapes. Once again we assume periodic boundary conditions and consider continuous-time solvers that predict $\nicefrac{\partial u_j}{\partial t}$ and use an ODE integration algorithm to advance $u_j$ in time. 
We again use a FV representation where the domain $\Omega$ is divided into $N$ grid cells $\Omega_j$ with volume $|\Omega_j|$ and the solution average in the $j$th grid cell is a scalar $u_j$.
We introduce bracket notation where $\langle \bm a \rangle \coloneqq \nicefrac{1}{N|\Omega|}\sum_{j=1}^N a_j |\Omega_j|$ denotes an average over the domain while the inner product notation $\langle \bm a | \bm b \rangle \coloneqq \sum_{j=1}^N a_j b_j |\Omega_j|$.

Suppose that the rate of change of $\bm u_j$ is given by
\begin{equation}\label{eq:blackbox}
 \frac{d \bm u_j}{d t}  = \bm N_j(\bm u_j)
\end{equation}
where $\bm N_j(\bm u_j) \in \mathbb{R}^N$ is an arbitrary update function. A machine learned solver would use ML to predict $\bm N_j$. Note that \cref{eq:blackbox} does not guarantee mass conservation by construction. Ensuring conservation of mass and the non-increasing $\ell_2$-norm property therefore requires modifying $\bm N_j$. Assuming periodic BCs, conservation of mass requires
\begin{equation}\label{eq:fv_arb_mass}
    \frac{d}{dt} \sum_{j=1}^N u_j |\Omega_j| = \sum_{j=1}^N \frac{d u_j}{d t} |\Omega_j| = \sum_{j=1}^N N_j  |\Omega_j|=0
\end{equation}
and non-increasing discrete $\ell_2$-norm requires
\begin{equation}\label{eq:fv_arb_l2}
    \frac{d}{dt} \sum_{j=1}^N \frac{1}{2} u_j^2 |\Omega_j| = \sum_{j=1}^N u_j \frac{d u_j}{d t} |\Omega_j| =  \sum_j u_j N_j |\Omega_j| \le 0.
\end{equation}
In vector-bracket notation, \cref{eq:fv_arb_mass,eq:fv_arb_l2} can be written as
$\langle \bm N_j \rangle = 0$
and
$\langle \bm u_j | \bm N_j \rangle \le 0$.
 These conditions will be satisfied if the following transformation is applied to $\bm N_j$:
\begin{equation}\label{eq:black_box_stability}
\begin{split}
    \bm U_j \coloneqq \bm u_j - &\langle \bm u_j\rangle \hspace{1.0cm} \bm M_j \coloneqq \bm N_j - \langle \bm N_j\rangle \hspace{1.0cm} \frac{d\ell_2^{\textnormal{old}}}{dt} = \langle \bm U_j | \bm M_j\rangle \\
     &\bm N_j \Rightarrow \bm M_j + \bigg(\frac{d\ell_2^{\textnormal{new}}}{dt} - \frac{d\ell_2^{\textnormal{old}}}{dt}\bigg)\frac{\bm G_j(\bm u_j)}{\langle \bm U_j | \bm G_j(\bm u_j)\rangle}
\end{split}
\end{equation}
for any $\nicefrac{d\ell_2^{\textnormal{new}}}{dt} \le 0$ and any finite function $\bm G_j(\bm u_j)$ where $\langle \bm G_j(\bm u_j) \rangle = 0$ and $\langle \bm G_j(\bm u_j) | \bm U_j \rangle \ne 0$. The choice $G_j(\bm u_j) = (\nabla^2 u)_j$ adds a spatially constant diffusion coefficient. In 1D with a spatially uniform grid, $(\nabla^2 u)_j = u_{j+1} - 2 u_j + u_{j-1}$.

\begin{figure}
  \centering
  \includegraphics[width=\textwidth]{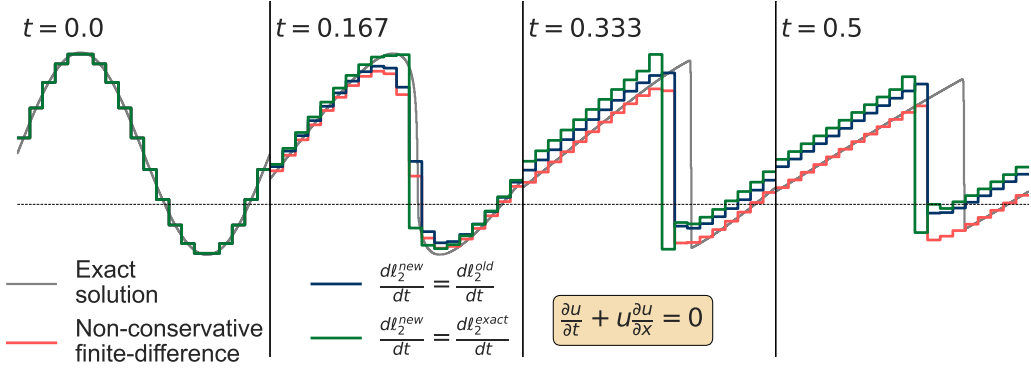}
  \caption{Using \cref{eq:black_box_stability}, we modify the time-derivative of a non-conservative finite-difference scheme. While the original finite-difference scheme (red) does not conserve the discrete mass $\sum_{j=1}^N u_j \Delta x$, the modified schemes (blue and green) conserve the discrete mass. All of these schemes decay the discrete $\ell_2$-norm. None of these finite-difference schemes result in a shock front traveling at the correct speed.} 
  \label{fig:burgers_nonconservative}
\end{figure}

We now demonstrate the effect of modifying a solver using \cref{eq:black_box_stability} to ensure the solver inherits discrete analogues of conservation of mass and non-increasing $\ell_2$-norm. We again solve the inviscid Burgers' equation $\frac{\partial u}{\partial t} + u \frac{\partial u}{\partial x } = 0$. The original scheme is the finite-difference scheme with
\begin{equation}
    N_j = -u_j (\delta u)_j\textnormal{, \hspace{0.5cm} } (\delta u)_j = \begin{cases}
        \nicefrac{u_j - u_{j-1}}{\Delta x}, & \text{if } u_j \ge 0\\
        \nicefrac{u_{j+1} - u_j}{\Delta x}, & \text{if } u_j < 0
        \end{cases}.
\end{equation}
We again use a 3rd-order SSPRK ODE integrator. The initial condition is $u_0(x) = 0.5 + \sin{x}$. The original finite-difference scheme, shown in red in \cref{fig:burgers_nonconservative}, does not conserve the discrete mass but preserves the non-increasing $\ell_2$-norm invariant. Using \cref{eq:black_box_stability} to modify $\bm N_j$, shown in blue and green in \cref{fig:burgers_nonconservative}, ensures that the discrete mass is conserved and maintains the property of non-increasing discrete $\ell_2$-norm. As $\nicefrac{d\ell_2^{\textnormal
old}}{dt}$ and $\nicefrac{d\ell_2^{\textnormal
exact}}{dt}$ are both negative and similar in magnitude, setting $\nicefrac{d\ell_2^{\textnormal
new}}{dt} = \nicefrac{d\ell_2^{\textnormal
exact}}{dt}$ (green) does not significantly improve accuracy over $\nicefrac{d\ell_2^{\textnormal
new}}{dt} = \nicefrac{d\ell_2^{\textnormal
old}}{dt}$. For all three of the schemes in \cref{fig:burgers_nonconservative}, the shock travels at the wrong speed.

\subsection{Discrete-Time FV Solvers with Arbitrary Time-Derivative \label{sec:discretetimeFV}}

\noindent We now show how to design invariant-preserving solvers when using the discrete-time update 
\begin{equation}\label{eq:discretetimeupdaterule}
    \bm u_{j}^{n+1} = \bm u_{j}^n + \Delta \bm u_j^n
\end{equation}
rather than the continuous-time update $\nicefrac{d \bm u_j}{d t} = \bm N_j(\bm u_j)$. The same error-correction strategy can be used to modify $\Delta \bm u_j^n$, though preserving the non-linear invariant requires solving a quadratic equation.
A discrete-time machine learned solver would use ML to predict $\Delta \bm u_j^n$. 
To ensure conservation of mass, we want
$\langle \bm u^{n+1}_j \rangle = \langle \bm u^{n}_j \rangle$
which requires that
\begin{equation}\label{eq:discretetime_masscondition}
    \langle \Delta \bm u^n_j \rangle = 0.
\end{equation}
To ensure that the discrete $\ell_2$-norm is non-increasing, we want
$\langle (\bm u_{j}^{n+1})^2 \rangle \le  \langle (\bm u_j^{n})^2 \rangle$. Suppose that we want the change of the discrete $\ell_2$ norm to be some scalar $\Delta \ell_2$ where $\Delta \ell_2 \le 0$. Thus, we want $\frac{1}{2}\langle (\bm u_{j}^{n+1})^2 \rangle  =  \frac{1}{2}\langle (\bm u_j^{n})^2 \rangle  + \Delta \ell_2$. Some simple algebra gives
\begin{equation}\label{eq:discretetime_l2condition}
     2 \langle \bm u_j^n |  \Delta \bm u_j^n \rangle +  \langle ( \Delta \bm u_j^n)^2 \rangle = 2 \Delta \ell_2 .
\end{equation}
\Cref{eq:discretetime_masscondition,eq:discretetime_l2condition} will be satisfied if the following transformation is made to $\Delta \bm u_j^n$:
\begin{equation}\label{eq:discretetime_updatedelta}
     \Delta \bm u_j^n \rightarrow \overline{\Delta \bm u_j^n}+ \epsilon \bm G_j(\bm u_j^n)
\end{equation}
where $\overline{\Delta \bm u_j^n} \coloneqq \bm{{\Delta \bm u}}_j^n - \langle \bm{{\Delta \bm u}}_j^n \rangle$, for any function $\bm G_j(\bm u_j)$ for which $\langle \bm G_j(\bm u_j^n) \rangle =0$ and $\epsilon$ is some yet-to-be-determined scalar. Plugging \cref{eq:discretetime_updatedelta} into \cref{eq:discretetime_l2condition} gives
\begin{equation}
    2\langle \bm u_j^n | \overline{\Delta \bm u_j^n} \rangle + \langle (\overline{\Delta \bm u_j^n})^2 \rangle + 2\epsilon\Big( \langle \bm u_j^n + \overline{\Delta \bm u_j^n} | \bm G_j \rangle \Big) + \epsilon^2 \langle  (\bm G_j)^2 \rangle = 2 \Delta \ell_2
\end{equation}
which is a quadratic equation for epsilon. \Cref{eq:discretetime_l2condition} will thus be satisfied if
\begin{equation}\label{eq:discretetime_epsilon}
    \epsilon = \frac{\langle \bm u_j^n + \overline{\Delta \bm u_j^n} | \bm G_j \rangle}{\langle  (\bm G_j)^2 \rangle}\Bigg[ -1 \pm \sqrt{1 - \frac{\langle  (\bm G_j)^2 \rangle \Big(2\langle \bm u_j^n | \overline{\Delta \bm u_j^n} \rangle + \langle (\overline{\Delta \bm u_j^n})^2 \rangle -  2 \Delta \ell_2\Big)}{\Big( \langle \bm u_j^n + \overline{\Delta \bm u_j^n} | \bm G_j \rangle \Big)^2}}\Bigg]
\end{equation}
To ensure that $\epsilon$ is small when ${{{\Delta \bm u}}_j^n}$ is small, we choose the plus sign for $\epsilon$ in \cref{eq:discretetime_epsilon}.
Modifying the discrete update $\Delta \bm u_j^n$ according to  \cref{eq:discretetime_updatedelta} where $\epsilon$ is given by the plus sign in \cref{eq:discretetime_epsilon} ensures that mass is conserved and the discrete $\ell_2$-norm changes by an amount $\Delta \ell_2$. 

Notice that \cref{eq:discretetime_epsilon} can have no solution. Depending on the values of $\bm u_j$, $\Delta \bm u_j^n$, and the hyperparameter $\bm G_j(\bm u_j)$, there is a minimum allowed value of $\Delta \ell_2$. We discuss the implications of this in \cref{sec:limitations}.

\begin{figure}
  \centering
  \includegraphics[width=\textwidth]{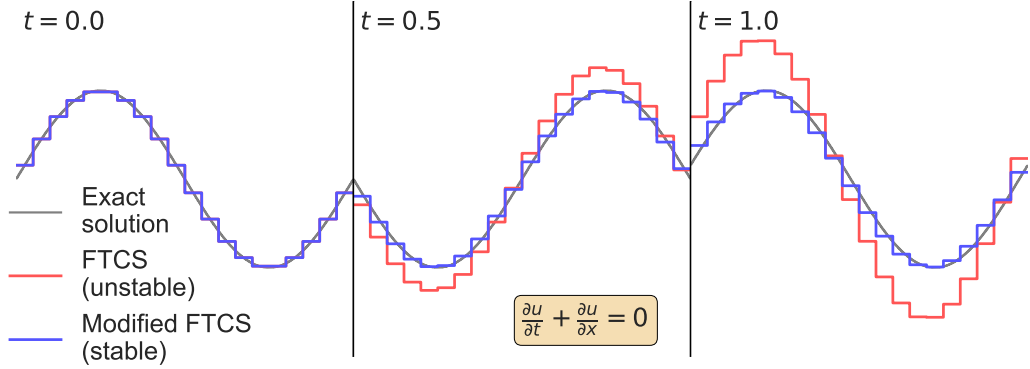}
  \caption{
    Solving the advection equation $\nicefrac{\partial u}{\partial t} + \nicefrac{\partial u}{\partial x} = 0$ with the forward-time centered-space (FTCS) update \cref{eq:ftcs} results in an $\ell_2$-norm increasing solution which is unconditionally unstable. Modifying the FTCS update using \cref{eq:discretetime_updatedelta,eq:discretetime_epsilon} allows us to control $\Delta \ell_2$, the change in the discrete $\ell_2$ norm. Setting $\Delta \ell_2 = 0$ (blue) results in an $\ell_2$-norm conserving and numerically stable solver.
  } 
  \label{fig:advection_ftcs}
\end{figure}

We now demonstrate the effect of using \cref{eq:discretetime_updatedelta,eq:discretetime_epsilon} to ensure that a discrete-time solver conserves mass and does not increase the discrete $\ell_2$-norm. We solve the advection equation $\frac{\partial u}{\partial t} + c \frac{\partial u}{\partial x} = 0$ with the forward-time, centered-space (FTCS) update 
\begin{equation}\label{eq:ftcs}
    u_{j}^{n+1} = u_{j}^n - \frac{ c\Delta t}{2\Delta x} (u_{j+1}^n - u_{j-1}^n)
\end{equation}
where $n$ is an index representing the values at the $n$th timestep. Our initial condition is $u_0(x) = \sin(x)$. We set $c=1$ and use a CFL number of $0.5$. We use periodic BCs in a domain of width $1$.
The exact solution to the advection equation is $u(x,t) = u_0(x - ct)$, which simply means that the solution translates to the right with speed $c$. The exact solution is shown in grey in \cref{fig:advection_ftcs}. 
FTCS conserves the discrete mass but increases the discrete $\ell_2$-norm for any $\Delta t$ and $\Delta x$ and is thus unconditionally unstable. The unstable FTCS solution is shown in red in \cref{fig:advection_ftcs}.
We can modify the FTCS update to control the change in the discrete $\ell_2$ norm using \cref{eq:discretetime_updatedelta,eq:discretetime_epsilon} with $G_j = u_{j+1} - 2 u_j + u_{j-1}$. For the 1D advection equation with smooth initial conditions we can set $\Delta \ell_2 = 0$, though for other PDEs we might want $\Delta \ell_2 < 0$. The $\ell_2$-norm conserving modified FTCS update is shown in blue in \cref{fig:advection_ftcs}. The modified update adds numerical diffusion to the FTCS update and results in a stable solver.

\subsection{Discontinuous Galerkin Solvers \label{sec:DG}}

\noindent In \cref{sec:continuoustimeFV}, we developed a technique for preserving invariants of any solver for \cref{eq:scalarhyperbolic} which uses a FV-like solution representation in the continuous-time limit. We now show how to do the same for any solver which uses a discontinuous galerkin (DG) solution representation. A machine learned DG solver would use ML to predict the time-derivative of the DG coefficients. We now give a brief introduction to DG methods \cite{dale_durran}.

\noindent \textbf{Discontinuous Galerkin}: Intuitively, the main difference between FV methods and DG methods is the solution representation: FV methods represent the solution as piecewise constant within each cell, while DG methods represent the solution as a polynomial within each cell.
With DG, we partition our domain $\Omega \in \mathbb{R}^d$ into $N$ cells $\{I_j\}_{j=1}^N$. In 1D, the solution representation within cell $I_j \in [x_{j-\nicefrac{1}{2}}, x_{j+\nicefrac{1}{2}}]$ is
\begin{equation}
    u_{j}(x,t) = \sum_{k=0}^p a_{jk}(t) \psi_k(x)
\end{equation}
where $\psi_k(x)$ are $p+1$ polynomial basis functions and $a_{jk}(t)$ are time-dependent coefficients. Notice that $u_j(x,t)$ is continuous within a cell but discontinuous across cell boundaries. The equation for the time-evolution of $a_{jk}(t)$ is found by minimizing the $\ell_2$-norm of the PDE residual within the subspace of basis functions spanned by $\psi_k(x)$. For \cref{eq:scalarhyperbolic}, this residual $E_j$ is given by
\begin{equation}
E_j = \int_{I_j}\bigg(\sum_{k=0}^p \dot{a}_{jk} \psi_k(x) + \frac{\partial f(u_j)}{\partial x}\bigg)^2\mathop{dx}.   
\end{equation}
The minimum of $E_j$ can be found from
\begin{equation}\label{eq:dg_weak_form}
    0 = \frac{\partial E_j}{\partial \dot{a}_{jk}} = 2 \int_{I_j} \psi_k \bigg(\sum_{k'=0}^p \dot{a}_{jk'} \psi_{k'} + \frac{\partial f}{\partial x}\bigg) \mathop{dx}.
\end{equation}
For general $\psi_k(x)$, \cref{eq:dg_weak_form} is a matrix equation for $\dot{a}_{jk}$. Typically, the $p+1$ basis functions are chosen to span the vector space of polynomials of degree $p$ and are orthogonal polynomials such that
\begin{equation}
    \int_{I_j}\psi_k(x)\psi_{k'}(x) dx = \Delta x_j \langle \psi_k | \psi_k \rangle {\delta_{kk'}}
\end{equation}
where $\langle \psi_k | \psi_k \rangle$ is a scalar and $\delta_{kk'}$ is the kronecker delta. In 1D, Legendre polynomials are usually chosen as basis functions so that $\langle \psi_k | \psi_k \rangle = 1/(2k+1)$. If orthogonal polynomials are chosen, then \cref{eq:dg_weak_form} can be inverted to solve for $\dot{a}_{jk}$:
\begin{equation}\label{eq:dg_before_ibp}
    \dot{a}_{jk} = -\frac{1}{\Delta x_j \langle \psi_k | \psi_k \rangle } \int_{I_j} \psi_k(x) \frac{\partial f}{\partial x} \mathop{dx}.
\end{equation}
The final step is to integrate by parts, giving
\begin{equation}\label{eq:dg_time_evolution}
     \dot{a}_{jk} =\frac{1}{\Delta x_j \langle \psi_k | \psi_k \rangle }\Big(- f_{j+\frac{1}{2}} \psi_k^+ + f_{j-\frac{1}{2}} \psi_k^- + \int_{I_j} f \frac{\partial \psi_k}{\partial x} \mathop{dx}\Big).
\end{equation}
where $\psi_k^+$ and $\psi_k^-$ are the values of $\psi_k$ at the right ($+$) and left ($-$) cell boundaries respectively. \Cref{eq:dg_time_evolution} has two terms: a volume term $\int_{I_j} f \frac{\partial \psi_k}{\partial x}dx$, and a boundary term $-f_{j+\nicefrac{1}{2}} \psi_k^+ + f_{j-\frac{1}{2}} \psi_k^-$. Depending on the form of \cref{eq:scalarhyperbolic}, the volume term can either be computed analytically or with a high-order quadrature. Because the discrete solution is discontinuous, the boundary term $f_{j+\frac{1}{2}}$ cannot be computed exactly and requires reconstructing the flux $f$ at cell boundaries. Note that if $k=0$ and $\psi_0(x)$ is chosen to be the zeroth Legendre polynomial $P_0(x)=1$, then \cref{eq:dg_time_evolution} reduces to the FV time-evolution equation \ref{eq:fv_a}. Note also that \cref{eq:dg_time_evolution} can be written as
\begin{equation}\label{eq:dg_general_update}
    \frac{d\bm a_{jk}}{dt} = \frac{1}{\Delta x \langle \psi_k | \psi_k \rangle} \bm N_{jk}
\end{equation}
where $\bm a_{jk}$ is the vector representation of the $N \times (p+1)$ solution coefficients. 
A machine learned DG solver would use ML to predict $\bm N_{jk}$.

\noindent \textbf{Discrete Invariants:} DG schemes conserve a discrete analogue $\sum_j \int_{I_j} u_j(x,t)\mathop{dx}$ of the continuous invariant $\int_\Omega u \mathop{dx}$ by construction. In a 1D periodic system, assuming the DG basis functions $\psi_k(x)$ are given by Legendre polynomials $P_k(x)$, we can see this with a short proof: 
\begin{equation}
\begin{split}
    \frac{d}{dt} \sum_{j=1}^N \int_{I_j} u_j(x,t) \mathop{dx} = &\sum_{j=1}^N \sum_{k=0}^p \dot{a}_{jk} \int_{I_j} P_k(x) \mathop{dx} \\= &\sum_j \dot{a}_{j0} \Delta x_j = -\sum_j(f_{j+\frac{1}{2}} - f_{j-\frac{1}{2}}) = 0.
\end{split}
\end{equation}
DG schemes do not automatically preserve discrete analogues of any of the non-linear invariants of the continuous PDE \cref{eq:scalarhyperbolic}. We now show how to modify a DG solver for \cref{eq:scalarhyperbolic} to ensure that the discrete $\ell_2$-norm is non-increasing in the continuous-time limit. We want
\begin{equation}
    \frac{d}{dt} \sum_j \int_{I_j} u_j(x,t)^2 \mathop{dx} \le 0.
\end{equation}
Using orthogonality of Legendre polynomials,
\begin{equation}\label{eq:DG_dl2dt}
    \begin{split}
        \frac{1}{2} \frac{d}{dt} \sum_j \int_{I_j} u_j(x,t)^2 dx = \frac{1}{2} \frac{d}{dt} \sum_j \int_{I_j} (\sum_k a_{jk} \psi_k)(\sum_{k'} a_{jk'} \psi_k') \mathop{dx} \\= \frac{1}{2}\frac{d}{dt}\sum_j \sum_k a_{jk}^2(t) \langle \psi_k | \psi_k \rangle \Delta x_j = \sum_j \sum_k a_{jk} \dot{a}_{jk} \langle \psi_k | \psi_k \rangle \Delta x_j.
    \end{split}
\end{equation}
Using \cref{eq:dg_general_update},
\begin{equation}
        \frac{1}{2} \frac{d}{dt} \sum_j \int_{I_j} u_j(x,t)^2 dx = \sum_j \sum_k a_{jk} N_{jk} = \langle \bm a_{jk} | \bm N_{jk} \rangle.
\end{equation}
Let us now define 
\begin{equation}
    \frac{d\ell_2^{\textnormal{old}}}{dt} \coloneqq \langle \bm a_{jk} | \bm N_{jk} \rangle
\end{equation}
as the original rate of change of the discrete $\ell_2$-norm, and $\nicefrac{d\ell_2^{\textnormal{new}}}{dt}$ as the desired rate of change of the discrete $\ell_2$-norm. To ensure non-increasing $\ell_2$-norm, we want $\nicefrac{d\ell_2^{\textnormal{new}}}{dt} \le 0$. To modify the rate of change of the discrete $\ell_2$-norm, we can add numerical diffusion (or anti-diffusion) to the time-derivative.
Unlike FV methods, where diffusion can be written as the sum of fluxes $f_{j+\nicefrac{1}{2}} \propto (u_{j+1} - u_j)$, DG diffusion requires computing both boundary terms and volume terms. The details of DG diffusion can be found in \cite{dgdiffusion}. For our purposes, it is sufficient to know that the time-derivative of the DG coefficients due to a diffusion term $\nu \nabla^2 u$ with diffusion coefficient $\nu$ can be written as
\begin{equation}\label{eq:dg_diffusion_term}
    \frac{d \bm a_{jk}}{dt} = \frac{\nu}{\Delta x \langle \psi_k | \psi_k \rangle} \bm N_{jk}^{\textnormal{diffusion}}.
\end{equation}
Using \cref{eq:DG_dl2dt}, the rate of change of the discrete $\ell_2$-norm due to the diffusion term is
\begin{equation}
    \frac{1}{2} \frac{d}{dt} \sum_j \int_{I_j} u_j(x,t)^2 dx = \nu \sum_{j}\sum_k a_{jk} N_{jk}^{\textnormal{diffusion}} = \nu \langle \bm a_{jk} | \bm N_{jk}^{\textnormal{diffusion}} \rangle.
\end{equation}
If we set
\begin{equation}\label{eq:dg_nu}
    \nu = \bigg(\frac{d\ell_2^{\textnormal{new}}}{dt}  - \frac{d\ell_2^{\textnormal{old}}}{dt}  \bigg) \frac{1}{\langle \bm a_{jk} | \bm N_{jk}^{\textnormal{diffusion}} \rangle}
\end{equation}
and add the diffusion term \cref{eq:dg_diffusion_term} to the original update equation \cref{eq:dg_general_update}, then the total rate of change of the discrete $\ell_2$-norm will be $\nicefrac{d\ell_2^{\textnormal{new}}}{dt}$.
Although the addition of a diffusion term is not the most general way of updating $\nicefrac{d\ell_2}{dt}$, it is straightforward to compute and has a clear physical interpretation.

\begin{figure}
    \centering
    \begin{subfigure}[b]{\textwidth}
        \centering
        \includegraphics[width=\textwidth]{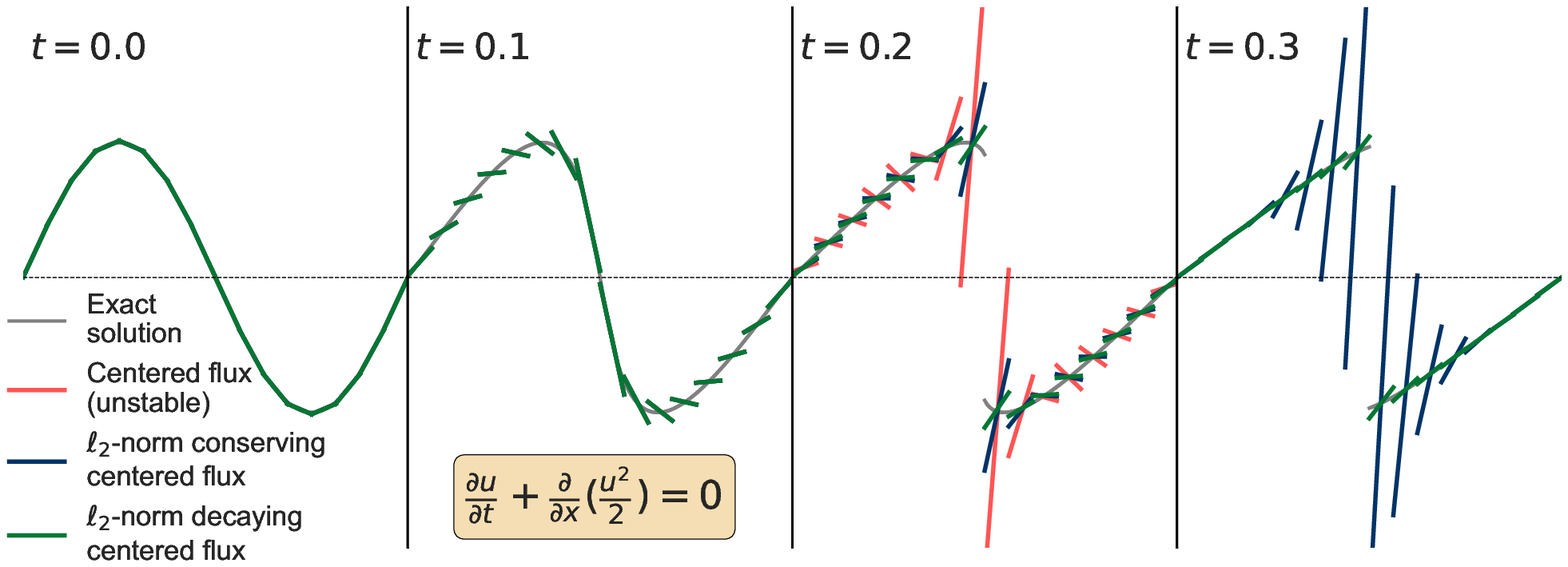}
        \caption{}
        \label{fig:dg_a}
    \end{subfigure}
    \begin{subfigure}[b]{\textwidth}
         \centering
         \includegraphics[width=\textwidth]{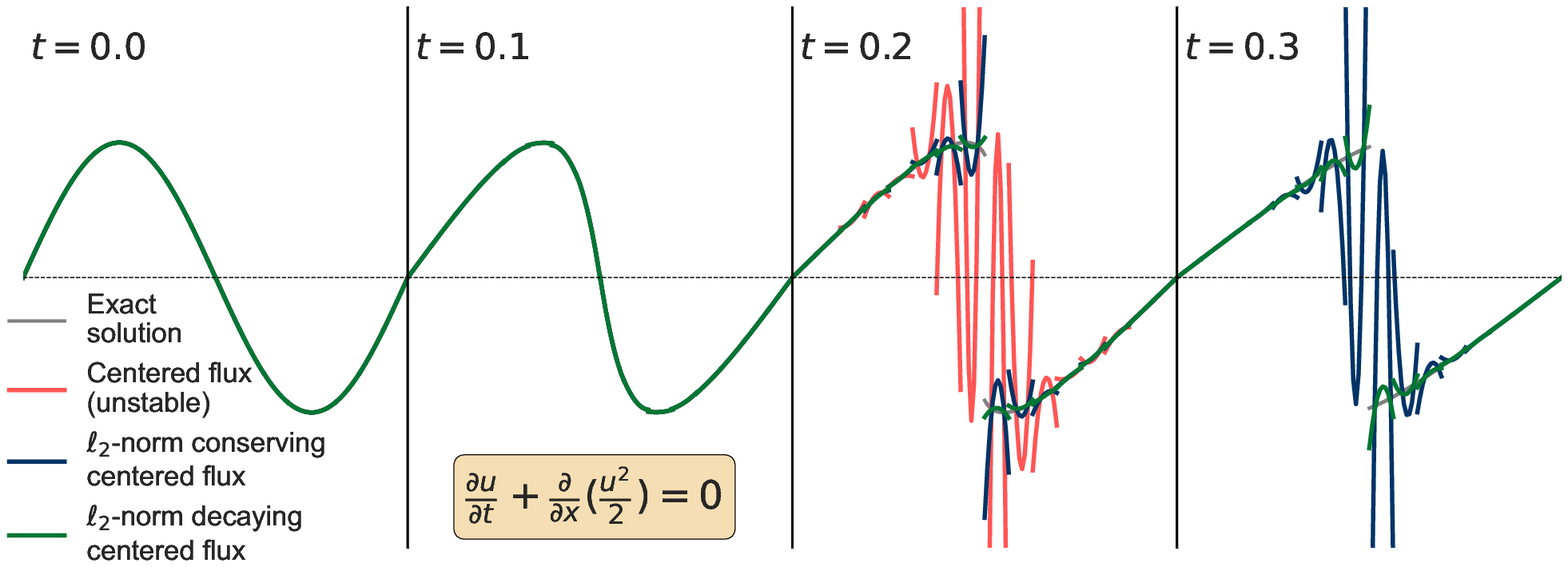}
         \caption{}
         \label{fig:dg_b}
     \end{subfigure}
    \caption{We add a diffusion term to two unstable DG solvers to control the rate of change of the discrete $\ell_2$-norm, thereby turning an unstable solver into a stable solver. Shown in red are standard DG solvers for the inviscid Burgers' equation with centered flux $f_{j+\nicefrac{1}{2}} = \frac{(u_{j+\nicefrac{1}{2}}^- + u_{j+\nicefrac{1}{2}}^+)^2}{8}$. This choice of flux increases the discrete $\ell_2$-norm, resulting in an unstable solver which blows up by $t=0.3$. (a) uses a degree-1 polynomial representation within each cell, while (b) uses a degree-2 polynomial representation within each cell. Adding a diffusion term with diffusion coefficient $\nu$ given by \cref{eq:dg_nu} controls the rate of change of the $\ell_2$-norm and results in a stable solver so long as $\nicefrac{d\ell_2^{\textnormal{new}}}{dt} \le 0$. Setting $\nicefrac{d\ell_2^{\textnormal{new}}}{dt} = 0$ (blue) gives large spurious oscillations and a stable but inaccurate solver, while setting $\nicefrac{d\ell_2^{\textnormal{new}}}{dt} = \nicefrac{d\ell_2^{\textnormal{exact}}}{dt}$ (green) adds additional diffusion and results in a more accurate solver.}
        \label{fig:dg_demo}
\end{figure}

Let us now demonstrate the effect of adding a diffusion term $\nu \nabla^2 u$ with $\nu$ given by \cref{eq:dg_nu} to an $\ell_2$-norm increasing DG solver. We use the same setup as in \cref{fig:fluxFV}. In \cref{fig:dg_demo} we show the time-evolution of DG methods with degree-1 polynomials ($p=1$, \cref{fig:dg_a}) and degree-2 polynomials ($p=2$, \cref{fig:dg_b}). The original DG solver uses the update \cref{eq:dg_time_evolution} with centered flux $f_{j+\nicefrac{1}{2}} = \frac{(u_{j+\nicefrac{1}{2}}^- + u_{j+\nicefrac{1}{2}}^+)^2}{8}$ where $u_{j+\nicefrac{1}{2}}^-$ is the solution just to the left of the ($j+\frac{1}{2}$)th cell boundary and $u_{j+\nicefrac{1}{2}}^+$ is the solution just to the right.
With this choice of flux, the original DG solvers are unstable and blow up by $t = 0.3$. By adding a diffusion term to the DG solver with diffusion coefficient given by \cref{eq:dg_nu}, we can control the rate of change of the $\ell_2$-norm and improve the accuracy of the solver. Like in \cref{fig:fluxFV}, setting $\nicefrac{d\ell_2^{\textnormal{new}}}{dt} = 0$ (blue) results in a highly oscillatory solution, but setting $\nicefrac{d\ell_2^{\textnormal{new}}}{dt} = \nicefrac{d\ell_2^{\textnormal{exact}}}{dt}$ damps many of the oscillations and results in a more accurate solution.


\subsection{Spectral Solvers \label{sec:fourier}}

\noindent In \cref{sec:fluxpredictingFV,sec:continuoustimeFV,sec:discretetimeFV} we considered solvers for \cref{sec:scalarhyperbolic} which represent the solution in a FV basis. In \cref{sec:DG}, we considered solvers which represent the solution in a DG basis. We now consider spectral solvers which represent the solution in a Fourier basis. In a 1D periodic domain where $x \in [0,L]$, the Fourier representation $\hat{u}$ of the solution $u$ is
\begin{equation}
    \hat{u}(x, t) = \sum_{m=-N}^{N} \tilde{u}_{m}(t) e^{\frac{2\pi i m x} {L}}.
\end{equation}
The coefficients $\tilde{u}_m = u_m^r + i u_m^i \in \mathbb{C}$. To ensure that $\hat{u}(x,t)\in \mathbb{R}$, we require $\tilde{u}_{-m} = \tilde{u}_{m}^*$ which gives $2N+1$ degrees of freedom in the solution representation. We consider the update equation
\begin{equation}
    \frac{d \tilde{\bm u}_m}{d t} = \tilde{\bm N}_m
\end{equation}
where $\tilde{\bm u}_m \in \mathbb{R}^{2N+1}$ is a vector representation of the $2N+1$ degrees of freedom in the complex solution coefficients.

\noindent \textbf{Discrete Invariants}: To ensure conservation of mass, we require that
\begin{equation}
    \frac{d}{dt} \int_0^L \hat{u}(x, t) \mathop{dx} = \int_0^L \sum_{m=-N}^N\frac{d\tilde{u}_m}{dt} e^{\frac{2\pi i m x} {L}} \mathop{dx} = L\tilde{N}_0 = 0.
\end{equation}
The rate of change of the $0$th Fourier coefficient must be zero. To ensure that the $\ell_2$-norm is non-increasing, we require that
\begin{equation}
    \frac{d}{dt}\int_0^L\frac{1}{2}|\hat{u}(x,t)|^2 dx 
    \le 0.
\end{equation}
We then use the Plancherel theorem
\begin{equation}
    \int_0^L |\hat{u}(x,t)|^2 dx = L\sum_{m=-N}^N |\tilde{u}_m|^2.
\end{equation}
Using $|\tilde{u}_m|^2=|\tilde{u}_{-m}|^2$ and $\nicefrac{d\tilde{u}_0}{dt}=0$, we have
\begin{equation}
   \frac{L}{2}\frac{d}{dt}\sum_{m=-N}^N |\tilde{u}_m|^2 = L\frac{d}{dt}\sum_{m=1}^N |\tilde{u}_m|^2 = 2L \sum_{m=1}^N u_m^r \frac{d{u}_m^r}{dt} + {u}_m^i \frac{d{u}_m^i}{dt} \le 0.
\end{equation}
In vector notation, this can be written as
\begin{equation}
    2L \langle \tilde{\bm u}_m | \tilde{\bm N}_m \rangle \le 0.
\end{equation}
These conditions will be satisfied if the following transformations are applied to $\tilde{\bm N}_m$:
\begin{equation}\label{eq:fourier_stability}
\begin{split}
    &\tilde{N}_0 \Rightarrow 0 \hspace{2.0cm} \frac{d\ell_2^{\textnormal{old}}}{dt} = 2L \langle \tilde{\bm{u}}_m | \tilde{\bm N}_m \rangle
    \\
    &\tilde{\bm{N}}_m \Rightarrow \tilde{\bm N}_m + \bigg(\frac{d\ell_2^{\textnormal{new}}}{dt} - \frac{d\ell_2^{\textnormal{old}}}{dt}\bigg)\frac{\bm G_m(\tilde{\bm u}_m)}{2L \langle \tilde{\bm u}_m | \bm G_m(\tilde{\bm u}_m)\rangle}
\end{split}
\end{equation}
for any $\nicefrac{d\ell_2^{\textnormal{new}}}{dt} \le 0$ and any finite function $\bm G_m(\tilde{\bm u}_m)$ where $G_0 = 0$ and $\langle \bm G_m(\tilde{\bm u}_m) | \tilde{\bm u}_m \rangle \ne 0$.

\subsection{2D Incompressible Euler Equations \label{sec:energyconservation}}

\noindent In \cref{sec:fluxpredictingFV,sec:continuoustimeFV,sec:discretetimeFV,sec:DG,sec:fourier}, we designed mass-conserving and $\ell_2$-norm non-increasing solvers for the generic scalar hyperbolic PDE, \cref{eq:scalarhyperbolic}. In this section, we design invariant-preserving solvers for a specific scalar hyperbolic PDE, the 2D incompressible Euler equations \cref{eq:euler}.
As we learned in \cref{sec:standardsolversincompressible}, the incompressible Euler equations exactly conserve the mass $\int \chi \mathop{dx}\mathop{dy}$, the energy $\frac{1}{2}\int \bm u^2 \mathop{dx}\mathop{dy}$, and the enstrophy $\int \chi^2 \mathop{dx}\mathop{dy}$. We will now design an error-correcting algorithm which preserves discrete analogues of these three invariants.

Suppose we want to solve \cref{eq:euler} on a 2D periodic rectangular domain. We divide the domain $\Omega$ into $N_x \times N_y$ cells with indices $i \in [1, \dots, N_x]$ and $j \in [1, \dots, N_y]$, with average vorticity $\chi_{ij}$, and with uniform grid spacing $\Delta x = \nicefrac{L_x}{N_x}$ and $\Delta y = \nicefrac{L_y}{N_y}$. Each cell has volume $|\Omega_{ij}| = \Delta x \Delta y$. We again use the notation $\mathop{dz} = \mathop{dx}\mathop{dy}$. Suppose also that the rate of change of $\bm \chi_{i,j}$ is given by
\begin{equation}\label{eq:update_incompressible}
    \frac{d\bm \chi_{i,j}}{dt} = \bm N_{i,j}(\bm \chi_{i,j})
\end{equation}
where $\bm N_{i,j} \in \mathbb{R}^{N_x \times N_y}$ is an arbitrary update function. A machine learned solver would use ML to predict $\bm N_{i,j}$ and would solve the elliptic equation $\nabla^2 \psi = - \chi$ using a standard FEM Poisson solver. Notice that \cref{eq:update_incompressible} is a subset of \cref{eq:blackbox} and that while $\bm \chi_{i,j}$ is represented in a discontinuous FV basis, $\bm \psi_{i,j}$ is represented in a continuous FEM basis.

\noindent \textbf{Discrete Invariants}: In \cref{sec:continuoustimeFV}, we showed that discrete conservation of mass requires $\langle \bm N_{i,j} \rangle = 0$ and non-increasing $\ell_2$-norm requires $\langle \bm u_{i,j} | \bm N_{i,j} \rangle \le 0$. The rate of change of the discrete energy can be computed using integration by parts and continuity of $\psi_{i,j}$ across cell boundaries:
\begin{equation}
\begin{split}
        \frac{d}{dt} \sum_{i,j} \int_{\Omega_{i,j}} \frac{1}{2} \bm \nabla \psi_{i,j} \cdot \bm \nabla \psi_{i,j} \mathop{dz} = 
    \sum_{i,j} \int_{\Omega_{i,j}}\bm \nabla \psi_{i,j} \cdot \bm \nabla \frac{\partial \psi_{i,j}}{\partial t} \mathop{dz}
    \\
    = -\sum_{i,j}\int_{\Omega_{i,j}} \psi_{i,j} \frac{\partial}{\partial t} \nabla^2 \psi_{i,j} + \int_{\partial \Omega_{i,j}} \psi_{i,j} \frac{\partial}{\partial t} \bm \nabla \psi_{i,j} \cdot \mathop{d\bm s} \\= \sum_{i,j} \int_{\Omega_{i,j}} \psi_{i,j} \frac{\partial \chi_{i,j}}{\partial t}\mathop{dz} = \sum_{i,j}  \frac{\partial \chi_{i,j}}{\partial t} \int_{\Omega_{i,j}} \psi_{i,j} \mathop{dz}.
\end{split}
\end{equation}
To ensure that the discrete energy is conserved, we require that
\begin{equation}
    \frac{d}{dt} \sum_{i,j} \int_{\Omega_{i,j}} \frac{1}{2} \bm \nabla \psi_{i,j} \cdot \bm \nabla \psi_{i,j} \mathop{dz} = 
    \sum_{i,j} \frac{\partial \chi_{i,j}}{\partial t} \int_{\Omega_{i,j}} \psi_{ij} \mathop{dz} = \langle \overline{\bm \psi}_{i,j} | \bm N_{i,j} \rangle = 0
\end{equation}
where $\overline{\psi}_{ij}$ is the average value of $\psi$ within cell $\Omega_{ij}$. Conservation of mass, conservation of energy, and the non-increasing $\ell_2$-norm property will therefore all be guaranteed for \cref{eq:euler} if the following transformation is applied to $\bm N_{ij}$:
\begin{align*}
    \bm U_{i,j} = \bm \chi_{i,j} - \langle \bm \chi_{i,j} \rangle && \bm M_{i,j} = \bm N_{i,j} - \langle \bm N_{i,j} \rangle && \bm{\Bar{\phi}}_{i,j} = \bm{\Bar{\psi}}_{i,j} - \langle \bm{\Bar{\psi}}_{i,j} \rangle \\
    \bm W_{i,j} = \bm U_{i,j} - \frac{\langle \bm U_{i,j} | \bm{\Bar{\phi}}_{i,j}  \rangle}{\langle \bm{\Bar{\phi}}_{i,j}  | \bm{\Bar{\phi}}_{i,j}  \rangle} \bm{\Bar{\phi}}_{i,j} && \bm P_{i,j} = \bm M_{i,j} - \frac{\langle \bm M_{i,j} | \bm{\Bar{\phi}}_{i,j}  \rangle}{\langle \bm{\Bar{\phi}}_{i,j}  | \bm{\Bar{\phi}} _{i,j} \rangle} \bm{\Bar{\phi}}_{i,j} && \frac{d\ell_2^{\textnormal{old}}}{dt} = \langle \bm W_{i,j} | \bm P_{i,j}  \rangle
\end{align*}
\begin{equation}\label{eq:eulerenergycons}
    \bm N_{i,j} \Rightarrow \bm P_{i,j} + \bigg(\frac{d\ell_2^{\textnormal{new}}}{dt} - \frac{d\ell_2^{\textnormal{old}}}{dt}\bigg)\frac{\bm G(\bm \chi_{i,j})}{\langle \bm W_{i,j} | \bm G( \bm \chi_{i,j}) \rangle}
\end{equation}
for any $\nicefrac{d\ell_2^{\textnormal{new}}}{dt}\le 0$ and any non-constant scalar function $\bm G_{i,j}(\bm \chi_{i,j})$ for which $\langle \bm G_{i,j}(\bm \chi_{i,j}) \rangle = 0$, $\langle \bm G_{i,j}(\bm \chi_{i,j}) | \bm {\bar{\psi}}_{i,j} \rangle = 0$ and $\langle \bm W_{i,j} | \bm G_{i,j}(\bm \chi_{i,j}) \rangle \ne 0$. The physically motivated choice 
\begin{equation}\label{eq:g_euler}
    \bm G_{i,j}(\bm \chi_{i,j}) = (\nabla^2 \bm W)_{i,j} - \frac{\langle (\nabla^2 \bm W)_{i,j} | \bm{\Bar{\phi}}_{i,j}  \rangle}{\langle \bm{\Bar{\phi}}_{i,j}  | \bm{\Bar{\phi}}_{i,j}  \rangle} \bm{\Bar{\phi}}_{i,j}.
\end{equation}
corresponds to the addition of a spatially constant diffusion coefficient, projected into an energy-conserving subspace.

We now illustrate the effect of modifying a stable standard solver for the 2D incompressible Euler equations using \cref{eq:2d_stability_a,eq:2d_stability_b,eq:eulerenergycons}. This standard solver is the second-order MUSCL scheme with monotonized central (MC) flux limiters \cite{sweby,muscl}. We use a linear finite element (FE) solver for the poisson equation \cite{fe_solver} and an SSP-RK3 ODE integrator \cite{ssprk}. The MUSCL scheme does not have any provable guarantees of $\ell_2$-norm conservation, but in 1D the MUSCL scheme is provably TVD \cite{osher1985convergence} and in practice the MUSCL scheme tends to decay the discrete $\ell_2$-norm as well as the discrete energy.

\begin{figure}
    \centering
    \begin{subfigure}[b]{0.98\textwidth}
        \centering
        \includegraphics[width=\textwidth]{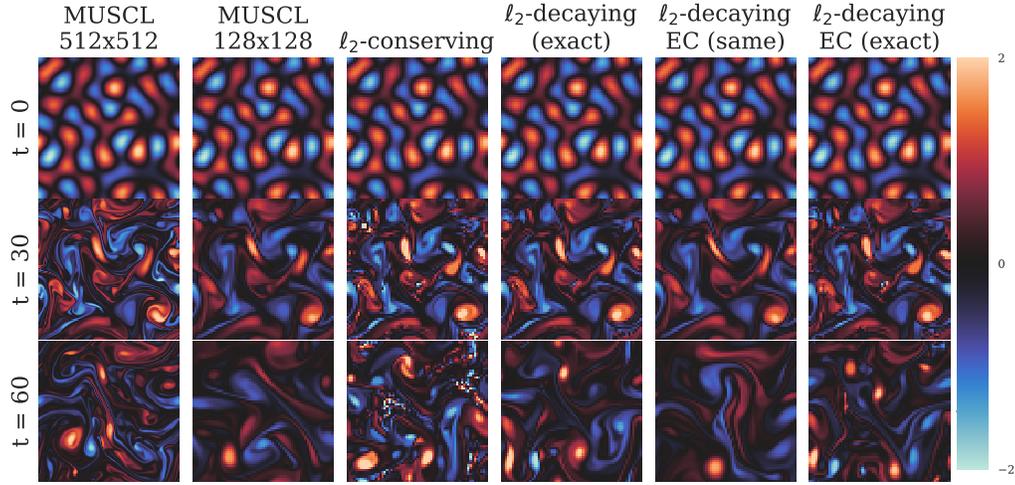}
        \caption{}
        \label{fig:euler_data}
    \end{subfigure}
    \begin{minipage}[c]{0.54\textwidth}
    \begin{subfigure}[b]{\textwidth}
         \centering
         \includegraphics[width=\textwidth]{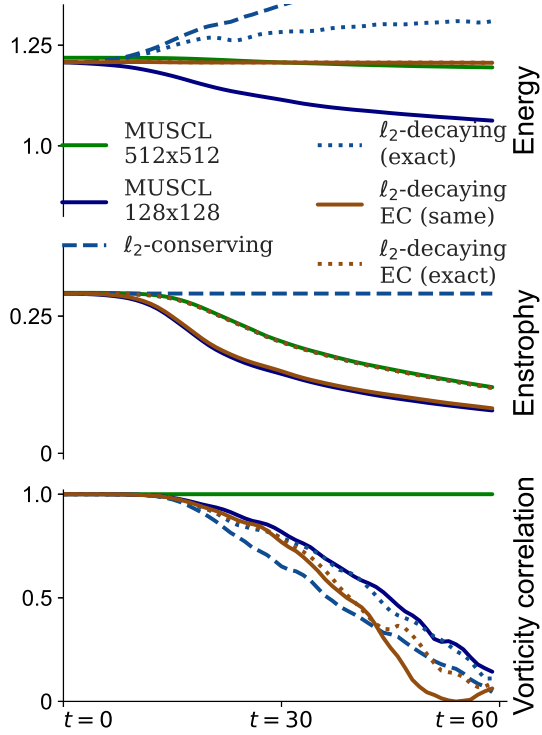}
         \caption{}
         \label{fig:euler_diag_demo}
     \end{subfigure}
     \end{minipage}
     \begin{minipage}[c]{0.44\textwidth}
    \caption{By modifying the time-derivative using \cref{eq:eulerenergycons}, we can design mass-conserving, energy-conserving, and enstrophy non-increasing solvers for the 2D incompressible Euler equations. (a) Images of the vorticity $\chi$ evolving under different numerical schemes for the incompressible Euler equations \cref{eq:euler}. The first and second columns show the standard MUSCL scheme at high and low resolution. The third and fourth columns show the low resolution MUSCL scheme, modified using \cref{eq:2d_stability_a,eq:2d_stability_b} to control the rate of change of the $\ell_2$-norm. The fifth and six columns use \cref{eq:eulerenergycons} to modify the low resolution MUSCL time-derivative to ensure energy conservation (EC) and again control the rate of change of the discrete $\ell_2$-norm. (b) Energy, enstrophy, and vorticity correlation over time. We use vorticity correlation as a benchmark measure of accuracy.}
    \label{fig:euler_2d}
    \end{minipage}
\end{figure}

In each of the six columns of \cref{fig:euler_data}, we see snapshots of the vorticity $\chi$ evolved using different numerical schemes. The first column is at high resolution ($512 \times 512$), while the other five columns are at low resolution ($128 \times 128$). The first and second columns use the unmodified MUSCL scheme. The third and fourth columns use \cref{eq:2d_stability_a,eq:2d_stability_b} to modify the MUSCL fluxes and set $\nicefrac{d\ell_2^{\textnormal{new}}}{dt}$. We use $\bm G_{i+\nicefrac{1}{2},j}^x = u_{i+1,j} - u_{ij}$ and $\bm G_{i,j+\nicefrac{1}{2}}^y = u_{i,j+1} - u_{ij}$. In the third column, we set $\nicefrac{d\ell_2^{\textnormal{new}}}{dt} = 0$ with $\nicefrac{d\ell_2^{\textnormal{new},x}}{dt} = \nicefrac{d\ell_2^{\textnormal{new},y}}{dt} = 0$. We find, similarly to \cref{fig:fluxFV}, that ensuring $\ell_2$-norm conservation introduces spurious high-$k$ oscillations. In the fourth column, we set 
$\nicefrac{d\ell_2^{\textnormal{new}}}{dt} = \nicefrac{d\ell_2^{\textnormal{exact}}}{dt}$, the rate of change of the discrete $\ell_2$-norm of the high resolution `exact' simulation.
We set $\nicefrac{d\ell_2^{\textnormal{new},x}}{dt} = \frac{1}{2} \nicefrac{d\ell_2^{\textnormal{exact}}}{dt}$, and $\nicefrac{d\ell_2^{\textnormal{new},y}}{dt} = \frac{1}{2}\nicefrac{d\ell_2^{\textnormal{exact}}}{dt}$. 
This allows for spurious oscillations to form, but much fewer than with $\nicefrac{d\ell_2^{\textnormal{new}}}{dt} = 0$. In the fifth and sixth columns, we modify the MUSCL time-derivative to enforce energy conservation (EC) using \cref{eq:eulerenergycons} with $\bm G_{i,j}$ set according to \cref{eq:g_euler}. In the fifth column we set $\nicefrac{d\ell_2^{\textnormal{new}}}{dt} =\nicefrac{d\ell_2^{\textnormal{old}}}{dt}$. In the sixth column we set $\nicefrac{d\ell_2^{\textnormal{new}}}{dt} = \nicefrac{d\ell_2^{\textnormal{exact}}}{dt}$. The energy-conserving schemes tend not to form spurious oscillations.

In \cref{fig:euler_diag_demo}, bottom row, we plot the vorticity correlation between the high resolution baseline and each of the five other schemes. Vorticity correlation has been used previously as a benchmark measure of accuracy for \cref{eq:euler} \cite{ml_accelerated_cfd}. We find that setting $\nicefrac{d\ell_2^{\textnormal{new}}}{dt} = 0$ worsens accuracy relative to the unmodified MUSCL scheme at the same resolution, while setting $\nicefrac{d\ell_2^{\textnormal{new}}}{dt} = \nicefrac{d\ell_2^{\textnormal{exact}}}{dt}$ neither helps nor harms accuracy. In \cref{fig:euler_diag_demo}, middle and top rows, we plot the discrete enstrophy $\frac{1}{2}\sum_{i,j}\int \chi_{i,j}^2 \Delta x \Delta y$ and discrete energy $\frac{1}{2}\sum_{i,j} \int (\bm u_{ij})^2\Delta x \Delta y$. The unmodified MUSCL schemes in the first and second columns of \cref{fig:euler_data} decay energy and enstrophy. 
The modified MUSCL schemes in the third and fourth columns of \cref{fig:euler_data} monotonically increase the discrete energy.
The energy-conserving schemes in the fifth and sixth columns of \cref{fig:euler_data} conserve the discrete energy.

\subsubsection{Enstrophy and Coarse Graining \label{sec:coarsegraining}}

\noindent The incompressible Euler equations in vorticity form exactly conserve the L2 norm of the solution, called the enstrophy. Naively, we might expect that a good PDE solver would conserve enstrophy as well. Yet algorithms that exactly conserve enstrophy on \cref{eq:euler}, such as the centered flux (not shown) or the modified MUSCL scheme with $\nicefrac{d\ell_2^{\textnormal{new}}}{dt} = 0$ (shown in \cref{fig:euler_2d}), perform significantly worse than schemes that allow enstrophy to decay.

To understand this puzzling result, we consider the relationship between the continuous enstrophy $\int \chi^2 \mathop{dx}\mathop{dy}$ and the discrete enstrophy $\sum_{i,j} \chi_{i,j}^2 \Delta x \Delta y$. Suppose that $\chi^{\textnormal{exact}} = \chi(x,y,t)$ is the exact solution to \cref{eq:euler}. As we know, $\chi^{\textnormal{exact}}$ has constant enstrophy. Now suppose that we coarse grain $\chi^\textnormal{exact}$, such that
$\chi^\textnormal{exact}_{i,j} = \int_{i,j} \chi^\textnormal{exact} \mathop{dx} \mathop{dy}$.
It turns out that, with very high probability, the discrete enstrophy of $\chi^{\textnormal{exact}}_{i,j}$ will decay in time \cite{basictypesofcoarsegraining}. This happens because $\chi^{\textnormal{exact}}$ tends to develop structures on a scale smaller than the grid size. These structures cannot be represented by $\chi_{i,j}$ and are replaced via coarse graining by a low-dimensional representation of the solution with lower enstrophy. 

Although the continuous equations for $\chi^\textnormal{exact}$ conserve enstrophy, the discrete equations for $\chi^\textnormal{exact}_{i,j}$ decay enstrophy.
Because machine learned PDE solvers solve discrete equations that are designed to approximate $\chi^{\textnormal{exact}}_{i,j}$, then machine learned PDE solvers should preserve the invariants of the discrete equations for $\chi^{\textnormal{exact}}_{i,j}$, not the invariants of the continuous equations for $\chi^{\textnormal{exact}}$.
As a result, solvers of \cref{eq:euler} should guarantee that enstrophy is non-increasing even though the continuous equations conserve enstrophy.

\section{Invariant-Preserving Error-Correcting Algorithms for Systems of Hyperbolic PDEs}
\label{sec:systemshyperbolic}

\noindent Unlike scalar hyperbolic PDEs, systems of hyperbolic PDEs are not guaranteed to all have the same non-linear invariants. While it is possible to design invariant-preserving algorithms for the generic scalar hyperbolic PDE \cref{eq:scalarhyperbolic}, invariant-preserving algorithms for systems of hyperbolic PDEs need to be tailored to the specific equation.

In this section, we design invariant-preserving error-correcting algorithms for an important example system: the compressible Euler equations of gas dynamics. We consider continuous-time, flux-predicting FV update rules in 1D. To produce invariant-preserving machine learned solvers, the procedure is the same as in \cref{sec:scalarhyperbolic}: at each timestep, if the invariants are not preserved apply an error-correcting algorithm to the update rule.

\subsection{Compressible Euler Equations\label{sec:compressible_euler}}

\noindent The 1D compressible Euler equations are given by \cref{eq:compressibleeuler}. They are in the form \cref{eq:system_conservation_form}. Recall from \cref{sec:standardsolvers} that some equations in the form \cref{eq:system_conservation_form} satisfy an entropy inequality \cref{eq:system_entropy_inequality} and that the compressible Euler equations satisfy an entropy inequality with generalized entropy function $\eta(s) = \rho g(s)$ with specific entropy $s = \log(\nicefrac{p}{\rho^\gamma})$ for which $\nicefrac{g''}{g'}\le \gamma^{-1}$. Recall that the choice $g(s) = e^{\nicefrac{s}{\gamma+1}}$ has the entropy variable $\bm w = (\nicefrac{\partial \eta}{\partial \bm u})^T$ given by \cref{eq:compressibleeuler_entropy_variable}.

Suppose we want to solve \cref{eq:compressibleeuler} for $\bm u = \begin{bmatrix}
    \rho, \rho v, E
\end{bmatrix}$
on a 1D domain with $x \in [0, L]$. We consider non-periodic BCs, though our results can easily be extended to periodic BCs. We use a FV discretization and divide the domain into $N$ cells of width $\Delta x = \nicefrac{L}{N}$ where the left and right boundaries of the $j$th cell for $j = 1, \dots, N$ are $x_{j-\nicefrac{1}{2}} = (j-1)\Delta x$ and $x_{j+\nicefrac{1}{2}}=j\Delta x$ respectively. We use a vector $\bm u_j = \begin{bmatrix}
        \rho_j,
        (\rho v)_j,
        E_j
    \end{bmatrix}$ to represent the solution average within each cell where $\bm u_j(t) \coloneqq \int_{x_{j-\nicefrac{1}{2}}}^{x_{j+\nicefrac{1}{2}}} \bm u(x, t)  \mathop{dx}$. We apply $\int_{x_{j-\nicefrac{1}{2}}}^{x_{j+\nicefrac{1}{2}}} (\dots) \mathop{dx}$ to \cref{eq:compressibleeuler} to derive the continuous-time FV update equation:
\begin{equation}\label{eq:compressible_euler_update}
    \frac{d \bm u_j}{d t} + \frac{1}{\Delta x} \bigg( \bm F_{j+\frac{1}{2}} - \bm F_{j-\frac{1}{2}}\bigg) = 0, \hspace{0.5cm} \bm F_{j+\frac{1}{2}} =  
    \begin{bmatrix}
        (\rho v)_{j+\frac{1}{2}} \\
        (\rho v^2 + p)_{j+\frac{1}{2}} \\
        (v(E + p))_{j+\frac{1}{2}}
    \end{bmatrix}.
\end{equation}
A machine learned solver would output the flux $\bm F_{j+\nicefrac{1}{2}}$ at $N-1$ cell boundaries for $j = 1, \dots, N-1$. $\bm F_{\nicefrac{1}{2}}$ and $\bm F_{N+\nicefrac{1}{2}}$ can be computed either from the boundary conditions or by using ML.

\noindent \textbf{Continuous Invariants}: As we saw in \cref{sec:standard_systems}, with non-periodic BCs the rate of change of the discrete mass is equal to the negative flux through the domain boundaries: $\frac{d}{dt}\int_{0}^L \bm u \mathop{dx} = \bm F(0) - \bm F(L)$. We also saw that the rate of change of the total entropy is greater than or equal to the entropy flux through the domain boundaries: $\frac{d}{dt} \int_{0}^L \eta(\bm u) \mathop{dx} \ge \psi(0) - \psi(L)$. The compressible Euler equations also maintain the positivity invariants $\rho \ge 0$ and $p = (\gamma - 1)(E - \frac{1}{2}\rho v^2) \ge 0$.

\textbf{Discrete Invariants}: With FV solvers, the rate of change of $\sum_{j=1}^N \bm u_j(t) \Delta x$ is equal to the negative flux through the boundaries $\bm F_{\nicefrac{1}{2}} - \bm F_{N+\nicefrac{1}{2}}$. Thus, the discrete analogue of the linear invariant $\int_\Omega u \mathop{dx}$ is preserved. We now show how to modify the predicted flux to ensure that the other two invariants are preserved. We want the discrete positivity invariants $\rho_j \ge 0$ and $p_j \ge 0$ to be maintained for all $j \in 1, \dots, N$. We also want the rate of change of a discrete analogue of the entropy to be greater than or equal to the entropy flux through the domain boundaries.

To ensure that the positivity invariants are preserved, we first limit $\bm F_{j+\nicefrac{1}{2}}$. One possible limiter introduced in \cite{Hu_2013} transforms $\bm F_{j+\nicefrac{1}{2}}$ for $j = 0, \dots, N$ according to the update
\begin{equation}\label{eq:compressible_euler_positivity_transformation}
    \bm F_{j+\nicefrac{1}{2}} \Rightarrow \theta_{j+\nicefrac{1}{2}} \bm F_{j+\nicefrac{1}{2}} + (1-\theta_{j+\nicefrac{1}{2}}) \bm F^{LF}_{j+\nicefrac{1}{2}}
\end{equation}
where $\bm F^{LF}_{j+\nicefrac{1}{2}}$ is the first-order Lax-Friedrichs flux and $0 \leq \theta_{j+\nicefrac{1}{2}} \leq 1$ is chosen to ensure positivity of $\rho_j$ and $p_j$. For details of how $\theta_{j+\nicefrac{1}{2}}$ is chosen, see \cite{Hu_2013}.

We now show how to further modify $\bm F_{j+\nicefrac{1}{2}}$ to ensure that the discrete entropy is greater than or equal to the entropy flux through the domain boundaries.
The discrete cell entropy $\eta_j = \rho_j g_j(s_j)$. We can choose $g_j(s_j) = e^{\nicefrac{s_j}{\gamma + 1}}$ with $s_j = \log{\nicefrac{p_j}{\rho_j^\gamma}}$ which results in the discrete entropy variable $\bm w_j = \frac{p^*_j}{p_j} \begin{bmatrix}
    E_j, -(\rho v)_j, \rho_j
\end{bmatrix}$ where $p^*_j = \frac{\gamma - 1}{\gamma + 1} (\nicefrac{p_j}{\rho_j^\gamma})^{\nicefrac{1}{\gamma + 1}}$ and $p_j = (\gamma - 1) (E_j - \nicefrac{(\rho v)_j^2}{2\rho_j} )$. Using $\nicefrac{d \eta_j}{d \bm u_j} = \bm w_j^T$, the rate of change of the discrete entropy is
\begin{equation}
    \frac{d}{dt}\sum_{j=1}^N \eta_j(\bm u_j) \Delta x = \sum_{j=1}^N \bm w_j^T \frac{d\bm u_j}{dt} \Delta x = -\sum_{j=1}^N \bm w_j^T (\bm F_{j+\frac{1}{2}} - \bm F_{j - \frac{1}{2}}).
\end{equation}
Using summation by parts, the rate of change of the discrete entropy is
\begin{equation}
        \frac{d}{dt} \sum_{j=1}^N \eta_j \Delta x = \bm F_{\frac{1}{2}}^T \bm w_1 - \bm F_{N+\frac{1}{2}}^T \bm w_N + \sum_{j=1}^{N-1} \bm F_{j+\frac{1}{2}}^T (\bm w_{j+1} - \bm w_j).
\end{equation}
To ensure that rate of change of the discrete entropy is greater than or equal to the entropy flux through the domain boundaries, we want
\begin{equation}\label{eq:compressibleeuler_entropyincreasing}
    \bm F_{\frac{1}{2}}^T \bm w_1 - \bm F_{N+\frac{1}{2}}^T \bm w_N + \sum_{j=1}^{N-1} \bm F_{j+\nicefrac{1}{2}}^T (\bm w_{j+1} - \bm w_j) \ge \psi(0) - \psi(L).
\end{equation}
In open systems, the entropy flux through the domain boundaries $\psi(0) - \psi(L)$ may or may not be known a priori. If it is not known a priori, we can make an estimate of the entropy flux using the boundary conditions.

Next, we transform $\bm F_{j+\nicefrac{1}{2}}$ to ensure that \cref{eq:compressibleeuler_entropyincreasing} is satisfied. We define scalars $\nicefrac{d\eta^{\textnormal{old}}}{dt} = \bm F_{\frac{1}{2}}^T \bm w_1 - \bm F_{N+\frac{1}{2}}^T \bm w_N + \sum_{j=1}^{N-1} \bm F_{j+\nicefrac{1}{2}}^T (\bm w_{j+1} - \bm w_j)$ and $\nicefrac{d\eta^{\textnormal{new}}}{dt} \ge \psi(0) - \psi(L)$. We then transform $\bm F_{j+\nicefrac{1}{2}}$ for $j = 1, \dots, N-1$ as follows:
\begin{equation}\label{eq:compressible_euler_entropy_transformation}
    \bm F_{j+\frac{1}{2}} \Rightarrow \bm F_{j+\frac{1}{2}} + \frac{ (\nicefrac{d\eta^{\textnormal{new}}}{dt} - \nicefrac{d\eta^{\textnormal{old}}}{dt})\bm G_{j+\nicefrac{1}{2}}}{\sum_{k=1}^{N-1} \bm G_{k+\nicefrac{1}{2}}^T(\bm w_{k+1} - \bm w_k)}
\end{equation}
for any finite non-constant function $\bm G_{j+\nicefrac{1}{2}}(\bm u)$ for which $\sum_{k=1}^{N-1} \bm G_{k+\nicefrac{1}{2}}^T(\bm w_{k+1} - \bm w_k) \ne 0$ and for which the addition of $\bm G_{j+\nicefrac{1}{2}}$ does not violate the positivity invariants.
The choice $\bm G_{j+\nicefrac{1}{2}}(\bm u) = \bm w_{j+1} - \bm w_j$ ensures that $\sum_{k=1}^{N-1} \bm G_{k+\nicefrac{1}{2}}^T(\bm w_{k+1} - \bm w_k) \ne 0$ so long as the discrete solution varies in space, but does not guarantee positivity of $\rho_j$ and $p_j$.
If we instead make the choice $\bm G_{j+\nicefrac{1}{2}}(\bm u) = \begin{bmatrix} 0, v_{j+1} -  v_j, p_{j+1} - p_j \end{bmatrix}^T$ corresponding to the physical diffusion terms $\nabla^2 u$ and $\nabla^2 p$ in the Navier-Stokes momentum and energy equations, we can guarantee positivity of $\rho_j$ and $p_j$ so long as $\nicefrac{d\eta^{\textnormal{new}}}{dt} > \nicefrac{d\eta^{\textnormal{old}}}{dt}$ (otherwise $\bm G_j$ is adding anti-diffusion) and a timestep restriction is satisfied. In our experiments, we empirically find that this choice also satisfies $\sum_{k=1}^{N-1} \bm G_{k+\nicefrac{1}{2}}^T(\bm w_{k+1} - \bm w_k) \ne 0$.

We now illustrate the effects of applying \cref{eq:compressible_euler_entropy_transformation} to an entropy-increasing, positivity-preserving numerical method. We use the Sod shock tube setup \cite{sod1978survey} with open (Dirichlet) BCs.
We apply the entropy-modifying algorithm \cref{eq:compressible_euler_entropy_transformation} to the MUSCL scheme \cite{muscl} with reconstruction in characteristic variables \cite{miyoshi2020short}.
We estimate the entropy flux through the boundaries ($\psi(0) - \psi(L)$) using the BCs and the formula $\psi = \rho g(s) v$. 
The MUSCL scheme is an entropy-increasing scheme, implying that $\nicefrac{\partial \eta^{\textnormal{old}}}{dt} > \psi(0) - \psi(L)$.
The entropy-modifying algorithm updates the rate of change of the entropy from $\nicefrac{\partial \eta^{\textnormal{old}}}{dt}$ to $\nicefrac{\partial \eta^{\textnormal{new}}}{dt}$, where $\nicefrac{\partial \eta^{\textnormal{new}}}{dt} = (\psi(0) - \psi(L)) + R (\nicefrac{\partial \eta^{\textnormal{old}}}{dt} - (\psi(0) - \psi(L)))$. We set $\bm G_{j+\nicefrac{1}{2}}(\bm u) = \begin{bmatrix} 0, v_{j+1} -  v_j, p_{j+1} - p_j \end{bmatrix}^T$.

\begin{figure}
    \centering
    \begin{subfigure}[b]{0.49\textwidth}
        \centering
        \includegraphics[width=0.9\textwidth]{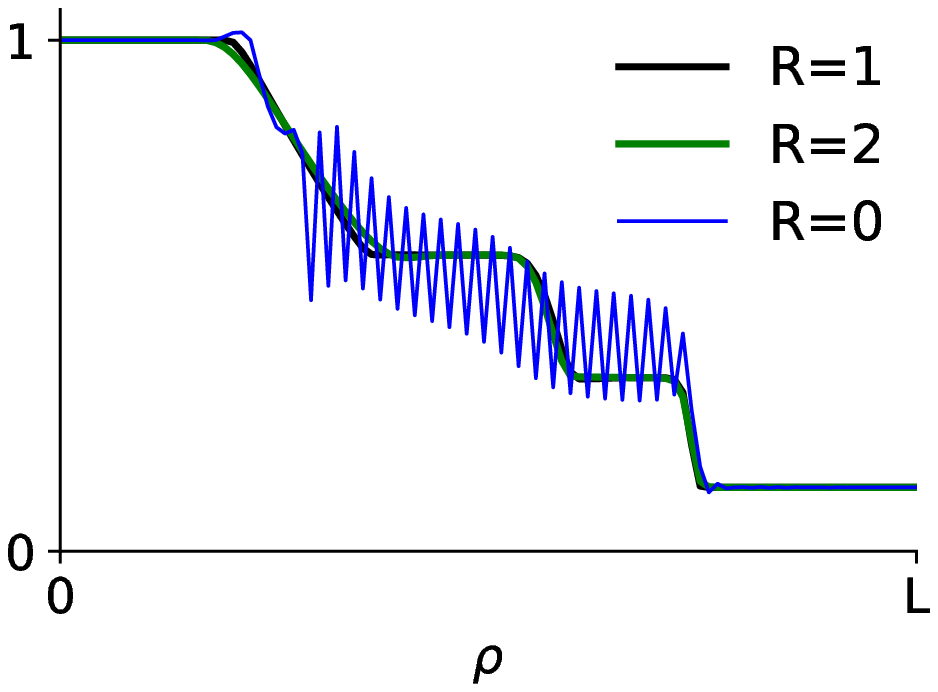}
        \caption{}
        \label{fig:3a}
    \end{subfigure}
    \begin{subfigure}[b]{0.49\textwidth}
         \centering
         \includegraphics[width=0.9\textwidth]{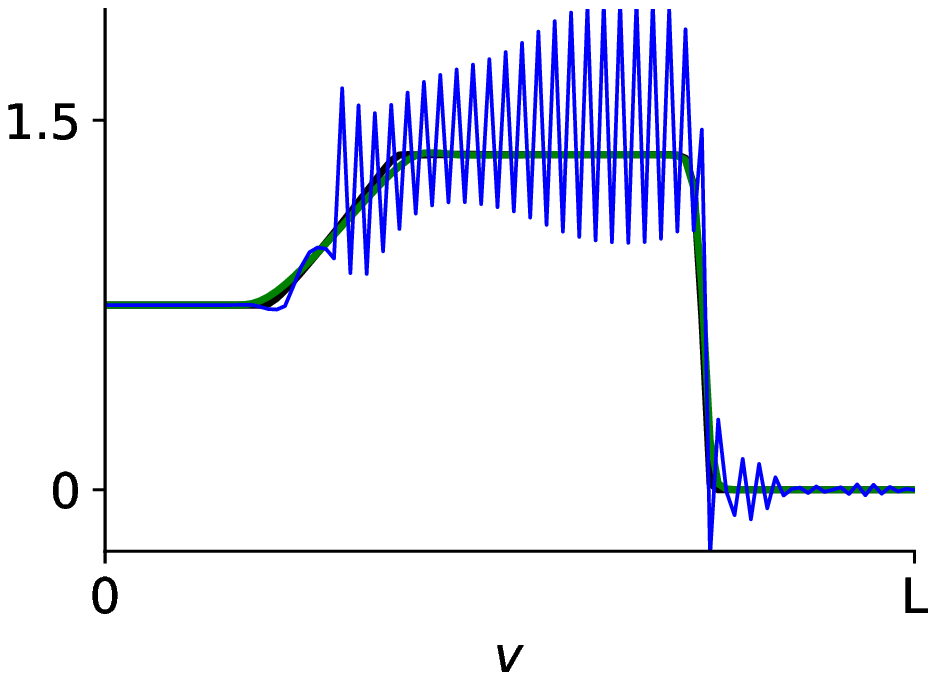}
         \caption{}
         \label{fig:3b}
     \end{subfigure}
    \begin{subfigure}[b]{0.49
    \textwidth}
         \centering
         \includegraphics[width=0.9\textwidth]{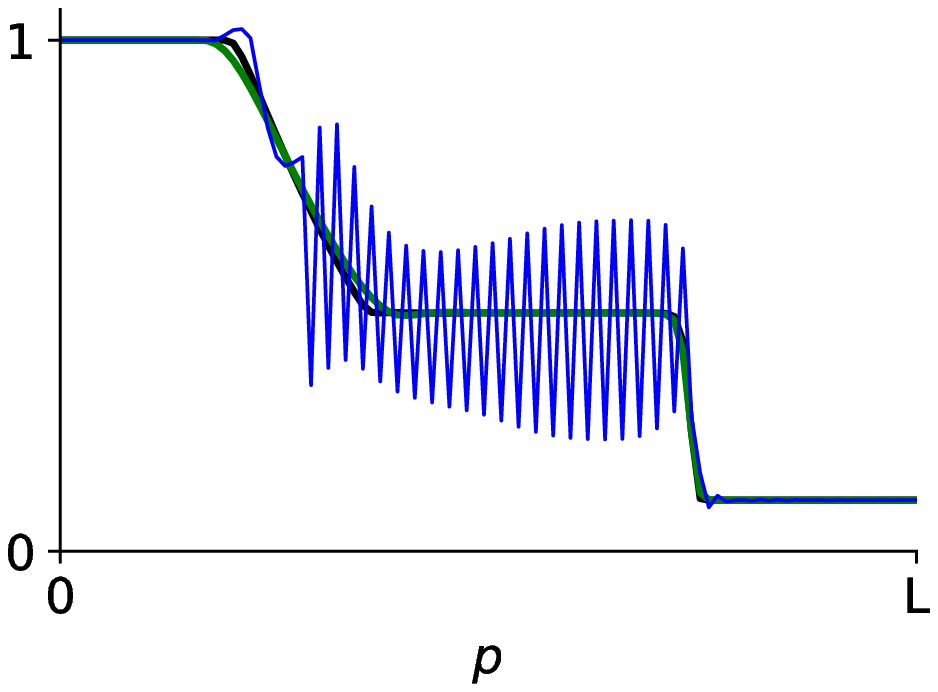}
         \caption{}
         \label{fig:3c}
     \end{subfigure}
    \caption{By adding diffusion and limiters to machine learned solvers using \cref{eq:compressible_euler_entropy_transformation,eq:compressible_euler_positivity_transformation} with $\nicefrac{\partial \eta^{\textnormal{new}}}{dt} > \psi(0) - \psi(L)$, we can design entropy-preserving and positivity-preserving solvers for the 1D compressible Euler equations. The discrete density $\rho$ in (a), velocity $v$ in (b), and pressure $p$ in (c) for three numerical methods for the 1D compressible Euler equations using a Sod shock tube setup with Dirichlet BCs. We modify the flux $\bm F_{j+\nicefrac{1}{2}}$ at cell boundaries of a standard MUSCL scheme using \cref{eq:compressible_euler_entropy_transformation}. We set $\nicefrac{\partial \eta^{\textnormal{new}}}{dt} = (\psi(0) - \psi(L)) + R (\nicefrac{\partial \eta^{\textnormal{old}}}{dt} - (\psi(0) - \psi(L)))$. The original MUSCL scheme ($R=1$) is shown in black. Increasing the rate of change of the discrete entropy ($R=2$), shown in green, adds diffusion to the original scheme. Decreasing the rate of change of the discrete entropy ($R=0$), shown in blue, adds anti-diffusion to the original scheme.}
        \label{fig:compressible_euler_demo}
\end{figure}

In \cref{fig:compressible_euler_demo} we plot the density $\rho$, velocity $v$, and pressure $p$ for three values of $R$. The original scheme ($R=1$) is plotted in black in \cref{fig:compressible_euler_demo}. Plotted in green is a scheme which increases entropy at a faster rate than the original scheme ($R=2$). Plotted in blue $(R=0)$ is a scheme with an entropy that changes only due to the entropy flux through the domain boundaries.
Notice that the scheme in blue adds anti-diffusion and is not guaranteed to preserve positivity, as $\nicefrac{d\eta^{\textnormal{new}}}{dt} < \nicefrac{d\eta^{\textnormal{old}}}{dt}$.
We do not show the result of setting $\nicefrac{d\eta^{\textnormal{new}}}{dt} < \nicefrac{d\eta^{\textnormal{BC}}}{dt}$ as the spurious oscillations result in $\rho_j < 0$ leading to NaNs in the reconstruction of the characteristic variables.

In summary: applying the transformation \cref{eq:compressible_euler_positivity_transformation} followed by the transformation \cref{eq:compressible_euler_entropy_transformation} with  $\nicefrac{d\eta^{\textnormal{new}}}{dt} \ge \psi(0) - \psi(L)$ to $\bm F_{j+\nicefrac{1}{2}}$ ensures that the discrete invariant $\sum_{j=1}^N \bm u_j(t) \Delta x$ is conserved, the positivity invariants $\rho_j\ge 0$ and $p_j \ge 0$ are maintained, and the discrete entropy $\sum_{j=1}^N \eta_j \Delta x$ is greater than or equal to the entropy flux through the domain boundaries.
These error-correcting transformations give us a procedure for designing invariant-preserving machine learned solvers. At each timestep or at each stage of an Runge-Kutta ODE integration, use \cref{eq:compressible_euler_positivity_transformation} to modify $\bm F_{j+\nicefrac{1}{2}}$ to ensure positivity of $\rho_j$ and $p_j$. If the resulting $\nicefrac{d\eta^{\textnormal{old}}}{dt} < \psi(0) - \psi(L)$, use \cref{eq:compressible_euler_entropy_transformation} with $\bm G_{j+\nicefrac{1}{2}}(\bm u) = \begin{bmatrix} 0, v_{j+1} -  v_j, p_{j+1} - p_j \end{bmatrix}^T$ to set $\nicefrac{d\eta^{\textnormal{new}}}{dt} = \psi(0) - \psi(L)$.
Because this transformation adds physically motivated diffusion terms to a positivity-preserving scheme, in the continuous-time limit the resulting machine learned solvers maintain positivity while ensuring that the discrete entropy is non-decreasing.

\section{Computational Verification}\label{sec:verification}

\noindent In \cref{sec:scalarhyperbolic,sec:systemshyperbolic} we algebraically derived invariant-preserving error-correcting algorithms. For most of these algorithms, we computationally verified that they do in fact preserve the correct invariants by applying them to standard numerical methods as shown in \cref{fig:fluxFV,fig:burgers_nonconservative,fig:advection_ftcs,fig:dg_demo,fig:euler_2d,fig:compressible_euler_demo}. In this section, we computationally verify two additional claims about these error-correcting algorithms.

First, in \cref{sec:why_standard_dont_work} we argue that, for certain invariants, standard approaches to preserving invariants will not work with machine learned solvers.
Intuitively, this is because standard approaches are either too restrictive or add too much numerical diffusion, degrading the accuracy of the solution.
Second, we claim that in closed or periodic systems the error-correcting algorithms we introduce do not degrade the accuracy of an already-accurate machine learned solver.
Intuitively, this is because an already-accurate solver will either satisfy or nearly satisfy the desired invariants, so the invariant-preserving correction will be small and only applied if necessary.
We verify this claim for the 1D advection equation in \cref{sec:why_standard_dont_work}, for the 1D Burgers' equation in \cref{sec:verification_burgers}, for the 2D incompressible Euler equation in \cref{sec:verification_2d_euler}, and for the 1D compressible Euler equations in \cref{sec:verification_compressible_euler}. We also demonstrate, using the 1D compressible Euler equations, how in an open system the accuracy of a machine learned solver can be degraded due to errors in the estimates of the rate of change of the invariant(s).

We emphasize that in this section, the purpose is not to design high-performing machine learned solvers but rather to use simple machine learned solvers trained to solve simple tasks as tools to illustrate the claims in the previous paragraph. Except for the 1D Burgers' solver in \cref{sec:verification_burgers}, we train extremely simple models. In particular, we don't use any ML-based strategies to improve performance, as the focus is not on the performance of the machine learned solvers. Instead, the focus is on the relative performance of the machine learned solvers with and without the invariant-preserving corrections.

\subsection{Why Standard Invariant-Preserving Methods Don't Work With ML\label{sec:why_standard_dont_work}}

Some invariants \textit{can} be preserved in machine learned solvers using standard approaches to invariant preservation without degrading accuracy. For example, the linear invariant  $\frac{d}{dt} \int_{\Omega} \bm u \mathop{dx} = 0$ can be preserved by predicting the flux across cell boundaries. Likewise, limiters can be used to ensure positivity of the solution. 
Yet for other invariants, standard approaches will not work with machine learned solvers. Why not? 

To be useful, machine learned PDE solvers must outperform standard numerical methods. To do so, invariant-preserving algorithms are needed which do not degrade accuracy at large $\Delta x$ and/or $\Delta t$. Standard approaches to preserving non-linear invariants are unstable at large $\Delta t$ due to the CFL condition and add local numerical diffusion proportional to $(\Delta x)^2$ \cite{artificial_viscosity}, thereby degrading accuracy at large $\Delta x$. As a result, machine learned solvers cannot simultaneously preserve non-linear invariants and outperform standard solvers using standard approaches.
The whole point of the invariant-preserving algorithms introduced in \cref{sec:scalarhyperbolic,sec:systemshyperbolic} is they guarantee invariant preservation at large $\Delta x$ and/or $\Delta t$ while adding the minimum correction needed to preserve the invariant.

We use an example to illustrate. We consider the simplest hyperbolic PDE, the 1D advection equation with periodic BCs.
As we will see, standard approaches can be used to design invariant-preserving machine learned PDE solvers, but those solvers cannot outperform standard methods. In contrast, the error-correcting algorithms we introduce preserve invariants without degrading accuracy. As a result, it becomes possible to design invariant-preserving machine learned solvers that outperform standard solvers.

\begin{figure}
  \centering
  \includegraphics[width=0.9\textwidth]{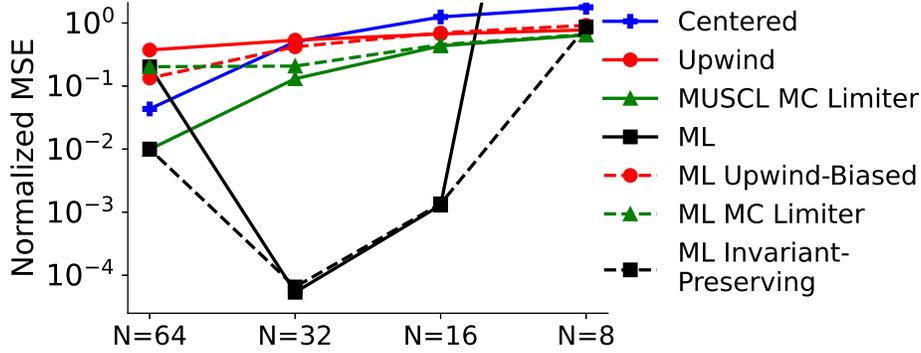}
  \caption{While there are various ways of designing invariant preserving machine learned PDE solvers, only the error-correction strategy we propose (black dotted line) preserves invariants without degrading the accuracy of an already-accurate machine learned solver (black line, $N=16$ and $N=32$ where $N$ is the number of spatial grid cells). Here we compare seven different methods for solving the 1D advection equation. Three are standard numerical methods (centered flux, upwind flux, and MUSCL flux) while four are machine learned solvers. After training the ML models with the same setup (see \ref{sec:appendix_verification_details}), we compute the normalized mean squared error (MSE) from $t=0$ to $t=1$ over 25 samples drawn from the training distribution.}
  \label{fig:compare_advection}
\end{figure}

We use the continuous-time FV update function \cref{eq:fv_a} and compare seven different choices for the flux at cell boundaries $f_{j+\nicefrac{1}{2}}$. Because they are flux-predicting methods, all seven solvers guarantee that the discrete mass is conserved. Three are standard numerical methods:
\begin{enumerate}
    \item The centered flux $f_{j+\nicefrac{1}{2}} = \frac{u_j + u_{j+1}}{2}$. This flux conserves the discrete $\ell_2$-norm.
    \item The upwind flux $f_{j+\nicefrac{1}{2}} = u_j$. This flux is TVD and decays the discrete $\ell_2$-norm. 
    \item The MUSCL flux with a Monotonized Central (MC) limiter. This flux is TVD.
\end{enumerate}
All three standard numerical methods (solvers 1, 2, 3) preserve one or more of the non-linear invariants and are thus numerically stable. The MUSCL scheme (solver 3) is the most accurate of the three standard numerical methods we consider. The other four methods are flux-predicting machine learned PDE solvers:
\begin{enumerate}
    \setcounter{enumi}{3}
    \item A machine learned solver which outputs $f_{j+\nicefrac{1}{2}}$. This solver is not guaranteed to conserve any non-linear invariants.
    \item An upwind-biased flux-predicting solver which outputs $\alpha_{j+\nicefrac{1}{2}} \ge 0$ where $f_{j+\nicefrac{1}{2}} = \alpha_{j+\nicefrac{1}{2}} f_{j+\nicefrac{1}{2}}^{\textnormal{Upwind}} + (1 - \alpha_{j+\nicefrac{1}{2}}) f_{j+\nicefrac{1}{2}}^{\textnormal{Centered}}$. This solver decays the $\ell_2$-norm.
    \item The same as solver 4, except with an MC flux limiter. This solver is TVD.
    \item The same as solver 4, but using the error-correcting algorithm \cref{eq:1d_stability} with $\nicefrac{d\ell_2^{\textnormal{new}}}{dt} \le 0$ to ensure that the discrete $\ell_2$-norm is non-increasing. 
\end{enumerate}
The upwind-biased solver (solver 5) and the flux-limited solver (solver 6) are examples of how standard approaches can be used to preserve invariants in machine learned PDE solvers. The error-corrected solver (solver 7) is an example of the strategy proposed in this paper. The details of the initial conditions, loss function, training data, and ML models are included in \ref{sec:appendix_verification_details}. In \cref{fig:compare_advection}, we plot the normalized mean squared error (MSE) for all seven solvers averaged over time from $t=0$ to $t=1$ averaged over 25 samples drawn from the training distribution. We compare solvers with $N$ grid cells, where $N=8$, 16, 32, and 64. 

The important takeaways from \cref{fig:compare_advection} are the following. The machine learned solver (solver 4) is highly accurate for $N=16$ and $N=32$. The upwind-biased machine learned solver (solver 5) decays the $\ell_2$-norm but is too constrained and too diffusive to outperform the standard solvers. The flux-limited machine learned solver (solver 6) is TVD but adds numerical diffusion proportional to $(\Delta x)^2$ to extremum and sharp gradients. At coarse resolution a high proportion of grid cells are either extremum or have sharp gradients. Thus flux-limited TVD-stable machine learned numerical methods operating at large $\Delta x$ will add large amounts of numerical diffusion to many of the grid cells which will degrade the accuracy of the solution.
In contrast, the invariant-preserving error-corrected machine learned solver (solver 7) adds the minimum amount of numerical diffusion necessary to ensure that $\nicefrac{d\ell_2^\textnormal{new}}{dt} \le 0$ and only does so when $\nicefrac{d\ell_2^\textnormal{old}}{dt} > 0$. Because the invariant-preserving correction adds the minimum amount of numerical diffusion necessary to preserve the invariant and only does so when the invariant is violated, the error-correcting algorithm preserves a discrete analogue of the non-increasing $\ell_2$-norm invariant without degrading the accuracy of an already-accurate solver ($N=16$ and $N=32$). 
For $N=8$ and $N=64$ the machine learned solver (solver 4) is unstable and increases the $\ell_2$-norm of the solution without bound, decreasing accuracy. There are a variety of possible causes of this poor performance, but the key takeaway is that the error-correcting solver (solver 7) improves the reliability and accuracy of a numerically unstable machine learned solver which has failed to preserve the desired invariants.

\subsection{Burgers' Equation \label{sec:verification_burgers}}

\noindent In this section, we verify for the 1D Burgers' equation that the error-correcting algorithms we introduce do not degrade the accuracy of an already-accurate machine learned solver. To do so, we need to train a machine learned solver to solve the 1D Burgers' equation.
Instead of designing our own solver, we attempt to replicate the solvers in fig 3c of \cite{data_driven_discretizations}; fig 3c compares the accuracy of a highly-accurate flux-predicting `data-driven discretization' 1D Burgers' solver with the accuracy of standard solvers. The solver, model, training, and evaluation are nearly identical to \cite{data_driven_discretizations}; details are included in \ref{sec:appendix_verification_details}. 

\begin{figure}
    \centering
    \includegraphics[width=0.86\textwidth]{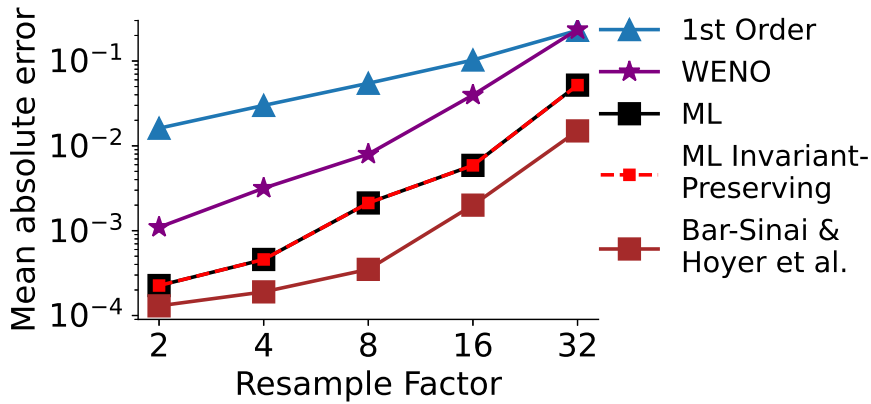}
    \caption{We demonstrate that the invariant-preserving error-correcting algorithm we introduce in \cref{sec:fluxpredictingFV} (red dotted line) doesn't degrade the accuracy of an already-accurate machine learned solver (black line).
    We plot the mean absolute error (MAE) averaged over 100 draws from the training distribution and over $t \le 15$.
    Instead of designing and training our own solver, we attempt to replicate the highly accurate 1D Burgers' machine learned solver in fig. 3 of \cite{data_driven_discretizations}.
    Our goal is for the accuracy of our replicated solver (black line) to roughly match the accuracy of the original solver (brown line). While we match the overall trend, our error is higher. Thus, our replication attempt is not fully successful.}
    \label{fig:ml_burgers}
\end{figure}

Our attempt to replicate fig 3c of \cite{data_driven_discretizations} is shown in \cref{fig:ml_burgers}. We plot the mean absolute error (MAE) for various solvers averaged over 100 samples drawn from the training distribution and over time $t$ less than 15. 
Our attempt at replicating the machine learned solver in \cite{data_driven_discretizations} is plotted in black.
In the dotted red line of \cref{fig:ml_burgers}, we apply the invariant-preserving error-correcting algorithm introduced in \cref{sec:fluxpredictingFV} to that machine learned solver. The accuracy of the invariant-preserving machine learned solver (red dotted line) is practically identical to that of the already-accurate machine learned solver (black line). This is because the already-accurate machine learned solver guarantees mass conservation and tends not to increase the discrete $\ell_2$-norm within its training distribution, so the correction is usually not needed. When the correction is needed, the error-correcting algorithm applies the minimum amount of numerical diffusion necessary to preserve the non-increasing $\ell_2$-norm invariant.

The MAE as originally reported in fig. 3 and fig. S8 of \cite{data_driven_discretizations} is plotted in brown. While are able to match the same trend as \cite{data_driven_discretizations}, we achieve worse accuracy. Unfortunately, our attempt at replicating fig. 3 of \cite{data_driven_discretizations} was only partially successful. We do not know the cause of the discrepancy. In any case, this failed replication attempt does suggest that it can be tricky to train machine learned solvers to achieve very high accuracy.

\subsection{2D Incompressible Euler Equation \label{sec:verification_2d_euler}}

\begin{figure}
    \centering
    \includegraphics[width=0.85\textwidth]{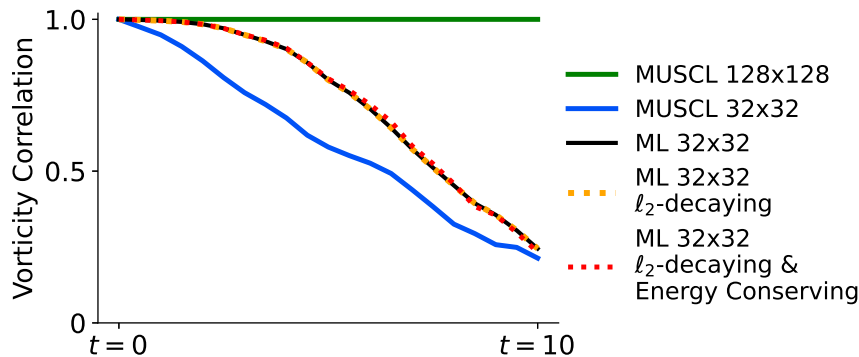}
    \caption{We demonstrate that the invariant-preserving error-correcting algorithms we introduce in \cref{sec:continuoustimeFV,sec:energyconservation} don't degrade the accuracy of an already-accurate machine learned solver (black line) for the 2D incompressible Euler equations. For each solver, we plot the correlation with the high-resolution `exact' solution; the correlation serves as a measure of accuracy. The low-resolution machine learned solver (black line) is more accurate than the low-resolution MUSCL scheme, but doesn't guarantee preservation of any of the invariants of the incompressible Euler equations. Applying error-correcting algorithms to either conserve mass and ensure non-increasing $\ell_2$-norm (\cref{eq:black_box_stability}, orange dotted line) or to conserve mass, conserve energy, and ensure non-increasing $\ell_2$-norm (\cref{eq:eulerenergycons}, red dotted line) do not degrade the accuracy of the machine learned solver. }
    \label{fig:corr_2d_euler}
\end{figure}

\noindent In this section, we verify for the 2D incompressible Euler equations that the invariant-preserving error-correcting algorithms we introduce do not degrade the accuracy of an already-accurate machine learned solver. 
The 2D incompressible Euler equations conserve mass, have non-increasing $\ell_2$-norm, and conserve energy. We consider the effect of two different invariant-preserving algorithms applied to a machine learned solver which is accurate at coarse resolution. First, we apply the mass conserving, $\ell_2$-norm non-increasing error-correcting algorithm in \cref{eq:black_box_stability} of \cref{sec:continuoustimeFV}. Second, we apply the mass conserving, energy-conserving, and $\ell_2$-norm non-increasing algorithm in \cref{eq:eulerenergycons} of \cref{sec:energyconservation}. We solve the 2D incompressible Euler equations with non-invariant forcing and diffusion terms, and use ML to approximate the invariant-preserving terms. We use standard methods to approximate the forcing and diffusion terms. The details of the solver, data generation, training and evaluation are in \ref{sec:appendix_verification_details}. \Cref{fig:corr_2d_euler} demonstrates that these invariant-preserving algorithms don't degrade the accuracy of an already-accurate machine learned solver.

\subsection{Compressible Euler Equations \label{sec:verification_compressible_euler}}

\noindent In this section, we verify for the 1D compressible Euler equations that in a periodic domain the invariant-preserving error-correcting algorithm we introduce in \cref{sec:compressible_euler} does not degrade the accuracy of an already-accurate machine learned solver, but that in a open system errors in the estimate of the rate of change of the invariant can degrade the accuracy of an already-accurate machine learned solver. The 1D compressible Euler equations conserve density, momentum, and energy, preserve positivity of density $\rho$ and pressure $p$,
and have non-decreasing entropy $\eta$. We train simple machine learned solvers in domains with periodic boundary conditions as well as domains with Dirichlet (open) boundary conditions. The initial conditions are given by a random draw from a (relatively simple) distribution of possible initial conditions; details of the data generation process, training procedure, and evaluation are given in \ref{sec:appendix_verification_details}.

\begin{figure}
    \centering
    \begin{subfigure}[b]{\textwidth}
        \centering
        \includegraphics[width=0.65\textwidth]{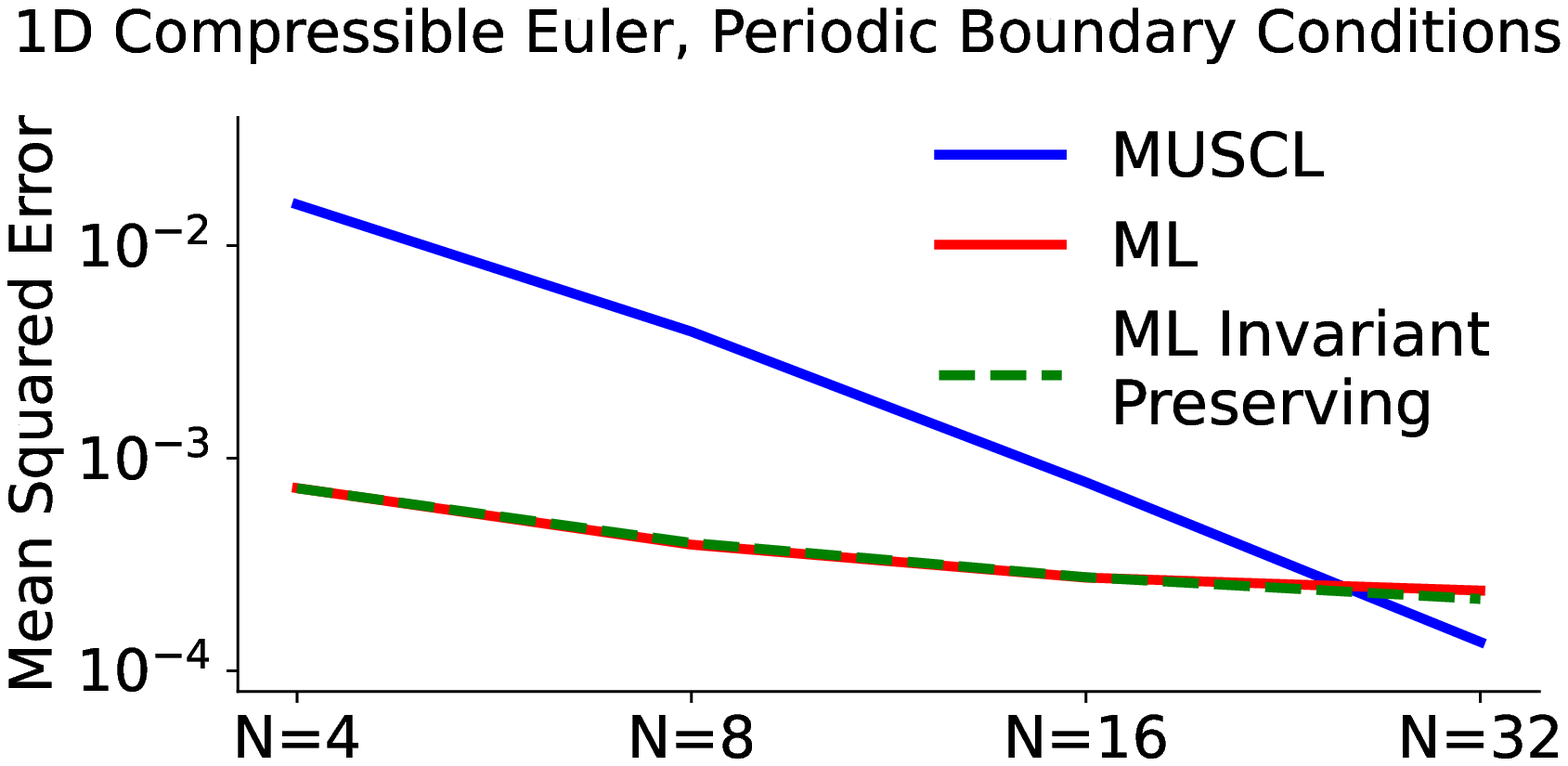}
        \caption{}
        \label{fig:1d_euler_periodic}
    \end{subfigure}
    \begin{subfigure}[b]{\textwidth}
         \centering
         \includegraphics[width=0.65\textwidth]{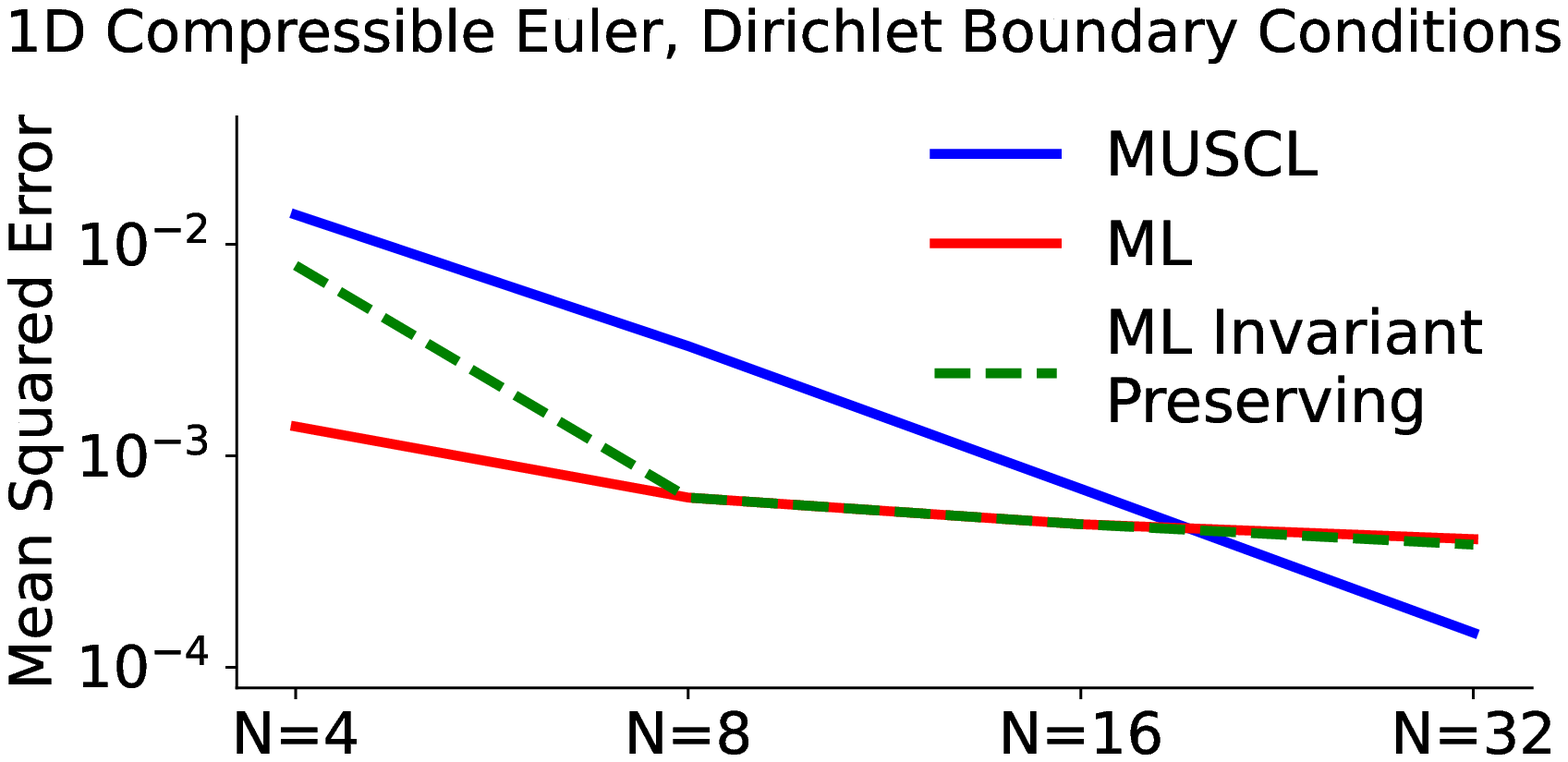}
         \caption{}
         \label{fig:1d_euler_open}
     \end{subfigure}
    \caption{We demonstrate that (a) the invariant-preserving algorithm for the 1D compressible Euler equations introduced in \cref{sec:compressible_euler} does not degrade the accuracy of a machine learned solver in periodic boundary conditions, and (b) that this algorithm can degrade the accuracy of a machine learned solver with open boundary conditions if the rate of change of the invariant is not estimated correctly. This error-correcting algorithm uses the transformation in \cref{eq:compressible_euler_positivity_transformation} to ensure positivity of $\rho_j$ and $p_j$, followed by the transformation in \cref{eq:compressible_euler_entropy_transformation} to ensure that entropy either is non-decreasing or is no less than the entropy flux through the domain boundaries.
    }
    \label{fig:1d_euler_verification}
\end{figure}

In \cref{fig:1d_euler_periodic}, we compare the performance of two machine learned solvers with a baseline MUSCL scheme (blue line) at various grid resolutions. The original machine learned solver (red line) learns a correction to the flux term of the MUSCL scheme. This solver conserves density, momentum, and energy, but doesn't guarantee positivity of $\rho_j$ or $p_j$ and doesn't guarantee that the discrete entropy will be non-decreasing. We apply the error-correcting algorithm introduced in \cref{eq:compressible_euler_positivity_transformation,eq:compressible_euler_entropy_transformation}, modified to support periodic BCs. This invariant-preserving machine learned solver (dotted green line) preserves the desired invariants without degrading the accuracy of the original solver.

In \cref{fig:1d_euler_open}, we do the same comparison except with Dirichlet boundary conditions. Due to the open boundary conditions, we do not have an exact estimate of the entropy flux through the boundaries. Thus, to use the invariant-preserving error-correcting algorithm \cref{eq:compressible_euler_positivity_transformation,eq:compressible_euler_entropy_transformation}, we need to estimate the entropy flux through the domain boundaries. Two possible estimates are to use the value from the Dirichlet boundary condition $\psi(0) = (\rho v)(0) g(s(0))$ or to use the value within the grid cell closest to the boundary $\psi(0) = (\rho v)_{1} g(s_1)$; we find that both estimates can add too much numerical diffusion and thereby degrade performance. We find better performance if we estimate the entropy flux to be the minimum of the two values, i.e., $\psi(0) = \textnormal{min}\{(\rho v)(0) g(s(0)), (\rho v)_{1} g(s_1)\}$. Applying that estimate to the error-correcting algorithm given by \cref{eq:compressible_euler_positivity_transformation}, in \cref{fig:1d_euler_open} we see that performance is degraded for the solver with $N=4$ grid cells but that the already-accurate machine learned solvers for $N=8$, $N=16$, and $N=32$ grid cells are able to guarantee positivity and satisfy the entropy-increasing invariant of \cref{eq:compressibleeuler_entropyincreasing} without degrading accuracy.

\section{Related Work}\label{sec:relatedwork}

\noindent \textbf{Invariant Preservation in Machine Learned Numerical Methods:} \cite{data_driven_discretizations,2d_learned_advection,ml_accelerated_cfd,stevens2020enhancement,stevens2020finitenet,stevens2022applications} use a finite volume representation of the solution to solve a variety of 1D and 2D PDEs. These so-called `hybrid' solvers predict the fluxes through cell boundaries and thus conserve mass by construction. They do not, however, preserve the non-increasing $\ell_2$-norm invariant and therefore do not guarantee stability. \cite{2d_learned_advection,ml_accelerated_cfd} attempt to promote stability by unrolling the loss function over multiple timesteps. \cite{ml_accelerated_cfd} trains a hybrid solver for the 2D incompressible Euler equations that ``remains stable during long simulations.'' This result is likely facilitated by the addition of physical diffusion to the PDE, which decays the $\ell_2$-norm at each timestep. \cite{thuerey_guaranteed_momentum} approximates a fluid with a system of particles and learns to predict the forces between particles using antisymmetric continuous convolutional layers; this ensures that the forces between particles are equal and opposite so that momentum is conserved. In such a system, we have found an error correction strategy could be used to conserve momentum \textit{and} energy. \cite{brandstetter2022clifford} solve Maxwell's equations in 3D and use Clifford Algebra and specially designed neural networks to ensure that geometric invariants relating the electric and magnetic fields are preserved. \cite{holloway2021acceleration} use ML to approximate the Boltzmann collision operator and apply the error-correcting algorithm introduced in \cite{zhang2018conservative}.
\cite{zhang2018conservative} introduce an error-correcting algorithm (intended for standard solvers) for the Boltzmann collision operator that enforces conservation of mass, momentum, and energy.
We have found an iterative error-correcting algorithm for the Boltzmann collision operator that preserves an additional invariant, non-decreasing entropy. We intend to discuss this algorithm in a future paper. \cite{richter2022neural} learn continuous divergence-free vector fields which satisfy incompressibility. \cite{alguacil2021predicting} applies an \textit{a posteriori} correction (i.e., error correction) to conserve energy in a deep learning-based solver for propagating acoustic waves; this correction results in improved performance.

\noindent \textbf{Other Machine Learned PDE Solvers:} \cite{learned_turbulence_modeling, solver-in-the-loop, ml_to_augment_sims, ml_accelerated_cfd} use convolutional neural networks to correct errors in low-resolution simulations; these hybrid solvers promote stability and improve accuracy by unrolling the loss function over multiple timesteps. \cite{deepmind_turbulence_sims} solves 2D and 3D hyperbolic PDEs using the `fully learned' update equation $u_{i,j}(t+\Delta t) = u_{i,j}(t) + N_{i,j}(u_{i,j}(t), \Delta t))$ where $N_{i,j}$ is a the output of a convolutional neural network; this update equation is identical to \cref{eq:discretetimeupdaterule}.
\cite{deepmind_turbulence_sims} attempts to ensure stability by adding noise to the training distribution and by using very large timesteps. 
\cite{message_passing_neural_pde_solvers} argues that instability in machine learned iterative numerical algorithms arises due to a \textit{distribution shift} where the distribution of training data differs from the outputs of the solver during inference due to small errors that accumulate over time. \cite{message_passing_neural_pde_solvers} uses the update equation $u_{i,j}(t+\ell\Delta t) = u_{i,j}(t) + \ell \Delta t N^\ell_{i,j}$ for $1 \le \ell \le K$ where $N^\ell_{i,j}$ is the output of a message passing graph neural network that predicts the next $K$ timesteps.
\cite{message_passing_neural_pde_solvers} attempts to ensure stability by modifying the loss function, adding random noise, and by predicting multiple timesteps into the future. A variety of papers have attempted to promote stability of dynamical systems that result from data-driven reduced order models, including by adding sparsity-promoting priors to a loss function \cite{Kaptanoglu_2021,erichson} and by constraining the eigenvalues of a learned Koopman operator \cite{koopman_stability}. 

\noindent \textbf{LES Models and Backscattering:} The objective of large eddy simulation (LES) is identical to that of many machine learned solvers: both attempt to find an accurate approximation to the solution of the PDE with fewer computational resources than classical numerical methods. Both also attempt to do so without resolving the smallest scales of the problem, relying on either an explicit or implicit subgrid model to do so  \cite{data_driven_discretizations,2d_learned_advection,ml_accelerated_cfd,deepmind_turbulence_sims,solver-in-the-loop,learned_turbulence_modeling,ml_to_augment_sims,Guan_2022,small_data_les, subgrid_modeling_2d_turbulence, les_closure, rose_yu_paper, ml_flux_limiters}.
Of particular relevance to the stability of subgrid models (both in LES and ML) are the concepts of `forward-scatter' and `backscatter'. In 2D LES turbulence, forward-scattering involves the transfer of enstrophy from resolved to unresolved scales, while backscattering involves the transfer of enstrophy from unresolved to resolved scales. Analysis across a wide range of flows demonstrates two important facts \cite{backscatter}. First, to be accurate a subgrid model must allow both forward-scatter and backscatter. This means that to be accurate a subgrid model must allow a discrete analogue of the entropy inequality \cref{eq:entropy_inequality} to be locally violated. Second, averaged over the entire domain there is always more forward-scatter than backscatter. If on average there were more backscatter than forward-scatter, then the subgrid model would be unstable \cite{Guan_2022}. In 2D turbulence, invariant-preserving error-correcting algorithms can be interpreted as ways of constraining a subgrid model to ensure that on average there is always at least as much forward-scatter as backscatter.

\section{Limitations \& Trade-offs \label{sec:limitations}}

\noindent Invariant preservation in machine learned solvers is, when applicable, free lunch. We know a priori our solution should preserve certain invariants, so by enforcing these invariants at each timestep we can improve the reliability of our solver without degrading accuracy. Nevertheless, there are some limitations and trade-offs to consider when designing and deploying these algorithms; we discuss eight.

First, the invariant-preserving updates in the error-correcting algorithms we introduce can all be derived analytically and applied in a single correction step. As a result, these algorithms are extremely efficient. However, in some situations it may be impossible to derive an analytic correction that preserves all the desired invariants. In such cases, iterative (and possibly less efficient) algorithms would need to be derived. In kinetic physics, for example, ensuring that the entropy $-\sum f \log{f}$ is non-decreasing seems to require an iterative error-correcting algorithm.

Second, many of the error-correcting algorithms we introduce add numerical diffusion in the continuous-time limit. If the machine learned solver is inaccurate and the algorithm needs to add large amounts of numerical diffusion, this can restrict the allowable timestep. Thus, inaccurate machine learned solvers might still be numerically unstable if the timestep is not modified to account for the additional diffusion.

Third, only one of the error-correcting algorithms we introduce (\cref{eq:discretetime_updatedelta,eq:discretetime_epsilon} in \cref{sec:discretetimeFV}) uses a discrete-time update rule. We have observed that for some invariants it seems to be more difficult (perhaps impossible) to derive analytic error-correcting algorithms with a discrete-time update than with a continuous-time update. 
Notice also that \cref{eq:discretetime_epsilon} has no solution if the change in the $\ell_2$-norm $\Delta \ell_2$ is below a solution- and solver-dependent minimum value. In practice, this means that inaccurate machine learned solvers that use a discrete-time update can still be numerically unstable. Specifically, since $\Delta \ell_2$ and $\Delta \hat{\bm u}$ both scale with the timestep $\Delta t$, more accurate machine learned solvers can guarantee numerical stability while using larger $\Delta t$, while less accurate solvers need to use smaller $\Delta t$ to do so.

Fourth, although designers of numerical algorithms usually would prefer to inherit discrete analogues of all of the properties of the continuous equation, it is usually impossible to design highly accurate numerical methods that inherit all of these properties.\footnote{Note also that the invariants of the continuous PDE are not necessarily the same as the invariants of the discrete equations. This was discussed in \cref{sec:coarsegraining}. }
This idea is formalized in Godunov's theorem \cite{godunov1959finite} for scalar hyperbolic PDEs, but is a pattern seen more generally in computational physics. For some PDEs there may be a trade-off between invariant preservation (i.e., reliability or robustness) and accuracy. Designers of numerical methods must choose which subset of invariants to preserve; being constrained by too many invariants can prevent a solver from making accurate predictions. In brief, it is possible to have \textit{too much} inductive bias.

Fifth, the invariant-preserving algorithms we propose are error-correcting algorithms which use global instead of local constraints. There is a trade-off associated with this decision: our algorithms violate the property that in hyperbolic PDEs information propagates at finite speed. While this trade-off may increase some readers conceptual discomfort with these algorithms, we believe it is a small price to pay for the benefits of numerical stability and improved reliability.

Sixth, while this error-correction strategy can still be applied to invariant-preserving PDEs containing non-invariant terms, the error-correcting algorithms should only be applied to the invariant-preserving terms. Thus, when solving non-invariant-preserving PDEs, machine learned solvers will need to separate out the contribution due to the invariant-preserving terms from the contribution due to the non-invariant-preserving terms.

Seventh, in open systems it can be difficult to estimate the fluxes through the boundary. Thus, while it is possible to design invariant-preserving error-correcting algorithms in open systems, if the rate of change of the invariant(s) is not estimated correctly then these algorithms can degrade accuracy. We suspect that this problem could be alleviated by using a higher density of grid cells near boundaries.

Eighth, while these error-correcting algorithms are designed to solve the problem of preserving invariants in machine learned solvers without degrading the accuracy of an already-accurate solver, they do not solve the problem of \textit{finding} accurate machine learned solvers. These algorithms adjust the update at each timestep if the solver violates invariants, but a solver which frequently commits large violations is likely to perform poorly. Alternatively, a solver could preserve the correct invariants but give inaccurate results. Building accurate, fast, and robust machine learned PDE solvers will require not only well-designed numerical methods but also well-engineered learning systems which consistently make accurate predictions about the time evolution of the solution.

\newpage

\appendix

\section{Details of Machine Learned Solvers\label{sec:appendix_verification_details}}

\noindent \textbf{Details for \cref{sec:why_standard_dont_work}}: 
We set $c=1$ and use periodic BCs. We use the continuous-time FV update equation \cref{eq:fv_a}. We use a SSPRK3 ODE integrator \cite{ssprk}. We choose the timestep $\Delta t$ using a CFL condition with a CFL number of 0.3. The initial conditions for both the training set and test set are draws from a sum-of-sines distribution
\begin{equation*}
    u_0(x) = \sum_{i=1}^{N^{\textnormal{modes}}} A_i \sin{\big(2\pi k_i x + \phi_i\big)}
\end{equation*}
where $N^\textnormal{modes} \sim \{1,2,3,4,5,6\}$ and $k_i \sim \{1,2,3,4\}$ are uniform draws from a set while $A_i \sim [-1.0,1.0]$ and $\phi_i \sim [0, 2\pi]$ are draws from uniform distributions. The loss function $L$ is given by computing the mean squared error (MSE) between the predicted time-derivative and the so-called `exact' time-derivative
\begin{equation*}
    L = \frac{1}{N_x} \sum_{j=1}^{N_x} \bigg (\frac{du_j(t)}{dt} - \frac{du_j^\textnormal{exact}(t)}{dt}\bigg)^2.
\end{equation*}
Both the `exact' solution $u_j^\textnormal{exact}(t)$ and the `exact' time-derivative $\nicefrac{du_j^\textnormal{exact}}{dt}$ are coarse-grained versions of a high-resolution simulation $u^\textnormal{exact}$, i.e., $u_j^\textnormal{exact}(t) = \int_{x_{j-\nicefrac{1}{2}}}^{x_{j+\nicefrac{1}{2}}} u^{\textnormal{exact}}(x,t) dx$. 
For each sample in the training data from the distribution of initial conditions, we take 50 snapshots evenly spaced in time from $t \in [0,1]$. We draw 100 samples, for a total of $5000$ snapshots in our training dataset. For each snapshot we store the exact trajectory $u_j^\textnormal{exact}$ and the exact time-derivative $\nicefrac{du_j^\textnormal{exact}}{dt}$.
Our machine learning models are periodic convolutional neural networks (CNNs) which are given the downsampled exact trajectories $u_j^{\textnormal{exact}}(t)$ as inputs. Solvers 4, 6, and 7 output the flux $f_{j+\nicefrac{1}{2}}$, while solver 5 outputs $\alpha_{j+\nicefrac{1}{2}}$. These outputs are then used to compute the predicted time-derivative $\nicefrac{du_j(t)}{dt}$.
All models train with a batch size of 32 and use the ADAM optimizer for 100 epochs over the training dataset with a learning rate of $1 \times 10^{-3}$, followed by another 100 epochs with a learning rate of $1 \times 10^{-4}$. Solver 7 uses $G_{j+\nicefrac{1}{2}}(\bm u_j) = u_{j+1} - u_j$.

\textbf{Details for \cref{sec:verification_burgers}}: Our goal is to solve the 1D Burgers' equation for $u(x,t) \in \mathbb{R}$ with diffusion and forcing:
\begin{equation}\label{eq:burgers_forcing}
    \frac{\partial u}{\partial t} + \frac{\partial}{\partial x} \bigg(\frac{u^2}{2}\bigg) = \nu \frac{\partial^2 u}{\partial x^2} + F(x, t).
\end{equation}
$x \in \Omega$ and $\Omega = [0, L]$. We set $\nu = 0.01$ and $L=2\pi$. We use periodic BCs. The initial conditions $u_0(x) = 0$. Each simulation in both the training and test data uses a randomized sum-of-sines forcing function
\begin{equation}
    F(x,t) = \sum_{m=1}^M A_m \sin{(2\pi k_m x/L - \omega_m t + \phi_m)}
\end{equation}
with $M=20$. The random variables $A_m \in [-0.5, 0.5]$, $\phi_m \in [0, 2\pi]$, and $\omega_m \in [-0.4, 0.4]$ are drawn from uniform distributions while $k_m$ is sampled uniformly from the set $\{3,4,5,6\}$. 

The machine learned solver and the standard solvers use the FV formulation. The domain $\Omega$ is divided into $N$ cells of width $\Delta x = \nicefrac{L}{N}$ and the solution average within each cell is represented by $u_j$ for $j = 1, \dots, N$. The FV update equations for \cref{eq:burgers_forcing} are
\begin{equation}
    \frac{\partial u_j}{\partial t} + \frac{\mathcal{J}_{j+\frac{1}{2}} -\mathcal{J}_{j-\frac{1}{2}} }{\Delta x} = F_j.
\end{equation}
The forcing 
\begin{equation}
    F_j(t) = \int_{x_j - \frac{\Delta x}{2}}^{x_j + \frac{\Delta x}{2}} F(x,t) \mathop{dx}
\end{equation}
is approximated using a 1-point 1st-order Guassian quadrature. For all solvers, the coefficients are advanced in time using a 3rd-order strong stability preserving Runge-Kutta \cite{ssprk} ODE integrator with a CFL factor of 0.3. Reducing the timestep further does not improve the accuracy of the baseline solvers.

The machine learned solver approximates the flux at each of the $N$ cell boundaries $\mathcal{J}_{j+\nicefrac{1}{2}}$ using the equation
\begin{equation}
    \mathcal{J}_{j+\frac{1}{2}} = \frac{1}{2}u_{j+\frac{1}{2}}^2 - \nu \bigg(\frac{\partial u}{\partial x}\bigg)_{j+\frac{1}{2}}.
\end{equation}
The reconstructed values of the solution $u_{j+\nicefrac{1}{2}}$ and its derivative $(\nicefrac{\partial u}{\partial x})_{j+\nicefrac{1}{2}}$ are approximated using the `data-driven discretization' approach introduced in \cite{data_driven_discretizations}. The data-driven discretization approach uses learned stencils $s^d_{j+\nicefrac{1}{2},k}$ to approximate the $d$th derivative of the solution:
\begin{equation}
    u_{j+\frac{1}{2}} = \sum_{k=1}^W s^0_{j+\frac{1}{2},k} u_{j-\frac{W}{2} + k} \quad\mathrm{and}\quad \bigg(\frac{\partial u}{\partial x}\bigg)_{j+\frac{1}{2}} = \sum_{k=1}^W s^1_{j+\frac{1}{2},k} u_{j-\frac{W}{2} + k}.
\end{equation}
We set the stencil width $W=6$. The learned stencils are computed as follows. First, a CNN maps the input array $u_j$ of length $N$ to two output arrays $\tilde{s}^0_{j+\nicefrac{1}{2},k}$ and $\tilde{s}^1_{j+\nicefrac{1}{2},k}$ of shape $(N, W)$. Second, for each of the $N$ stencil coefficients $\tilde{s}$, we project the length-$W$ vector of stencil coefficients into the null space of a matrix $M^d$, resulting in intermediate stencils $\bar{s}$. These projections can be written as 
\begin{equation}
    \bar{s}^0_{j+\frac{1}{2},k} = \tilde{s}^0_{j+\frac{1}{2},k} - \frac{1}{W}\sum_{l=1}^W \tilde{s}^0_{j+\frac{1}{2},l} \quad \mathrm{and} \quad \bar{s}^1_{j+\frac{1}{2},k} = \tilde{s}^1_{j+\frac{1}{2},k} - P^1_{k,l} \tilde{s}^1_{j+\frac{1}{2},l}
\end{equation}
where the projection matrix $P^1_{k,l}$ is given by
\begin{equation}
    P^1 = (M^{1})^T(M^1 (M^1)^T)^{-1} M^1
\end{equation}
where
\begin{equation}
    M^1 = \begin{bmatrix}
        1 & 1 & 1 & 1 & 1 & 1\\
        \frac{-5}{2} & \frac{-3}{2} & \frac{-1}{2} & \frac{1}{2} & \frac{3}{2} & \frac{5}{2}
    \end{bmatrix}
\end{equation}
Third, we add a `base' stencil $\hat{s}$ to the intermediate stencil $\bar{s}$. 
\begin{equation}
    {s}^0_{j+\frac{1}{2},k} = \bar{s}^0_{j+\frac{1}{2},k} + \hat{s}^0_k \quad \mathrm{and} \quad {s}^1_{j+\frac{1}{2},k} = \frac{\bar{s}^1_{j+\frac{1}{2},k} + \hat{s}^1_k}{\Delta x}
\end{equation}
where $\hat{s}^0_k = \begin{bmatrix}0 & 0 & \nicefrac{1}{2} & \nicefrac{1}{2} & 0 & 0 \end{bmatrix}$ and $\hat{s}^1_k = \begin{bmatrix}0 & 0 & -1 & 1 & 0 & 0 \end{bmatrix}$.
These three steps preserve formal 1st-order accuracy of the stencil coefficients. 

The CNN first applies two periodic convolutions with kernel size $K=5$, 32 channels, and ReLU activation function. The CNN then applies a periodic convolution with kernel size $K=4$ and $2 \times W$ output channels. A periodic convolution involves periodically padding the array with $K-1$ pixels then applying a `valid' convolution. The CNN has a receptive field of 12, half of which are on either side of the $j+\nicefrac{1}{2}$th cell boundary. The CNN weights are initialized using the LeCun normal initialization. The biases are initialized to zero.

The training data is generated using a high-resolution simulation with $N=512$ grid points and the WENO5 \cite{weno5} flux function. We generate data from 800 simulations, storing 10 snapshots per simulation taken every $t=0.5$ units of time, starting at $t=0$. We save the solution $u_j$ as well as the downsampled high-resolution `exact' time-derivative $\nicefrac{\partial u_j}{\partial t}$. We downsample by averaging $u_j$. 

Our loss function is given by the mean squared error between the downsampled high-resolution `exact' time-derivative and the time-derivative from the `data-driven discretization' time-derivative. We train with a batch size of $128$ times the downsampling factor. We use the Adam optimizer. We train for 20,000 steps with a learning rate of $3 \times 10^{-3}$ followed by 20,000 steps with a learning rate of $3 \times 10^{-4}$. 

\Cref{fig:ml_burgers} is computed by averaging the mean absolute error over 100 simulations in the test set. The `exact' solution is given by a high-resolution solver using the WENO flux. We initialize with $u_0(x)=0$, then compute the average error with each lower-resolution method from $t=0$ to $t=15$.

While our setup is nearly identical to the setup in \cite{data_driven_discretizations}, there are a few minor differences. First, \cite{data_driven_discretizations} incorrectly computes the Godunov flux, which leads to numerical instability of their baseline methods and worse performance. Second, we use a SSPRK3 ODE integrator while \cite{data_driven_discretizations} uses an order 3(2) adaptive time-stepping method from SciPy. Third, our CNN has a receptive field which is 12 points wide (symmetric around the $j+\nicefrac{1}{2}$th cell boundary) while the CNN in \cite{data_driven_discretizations} has a receptive field which is 13 points wide. We don't believe these minor differences are the root cause of the discrepancy in \cref{fig:ml_burgers} between our replication attempt (ML, black line) and figure 3 of \cite{data_driven_discretizations} (Bar-Sinai \& Hoyer et al., brown line). We do not know why we are unable to exactly replicate figure 3 of \cite{data_driven_discretizations}.

\noindent \textbf{Details for \cref{sec:verification_2d_euler}}: Our goal is to solve the 2D incompressible Euler's equations for $\chi(x,y,t)$ with forcing and diffusion:
\begin{align}
    \frac{\partial \chi}{\partial t} + \bm \nabla \cdot (\bm u \chi) = F(x,y,t) + \nu \nabla^2\chi \textnormal{,} && \bm u = \bm \nabla \psi \times \hat{e}_z\textnormal{,} && -\bm\nabla^2 \psi = \chi.
\end{align}
We choose periodic boundary conditions such that $x, y \in [0, L]$ and $L=2\pi$. We set $\nu = 10^{-3}$ and use the forcing function in \cite{ml_accelerated_cfd},
\begin{equation}
    F(x,y,t) = \frac{2 \pi k}{L} \cos{\bigg(\frac{2\pi k y}{L}\bigg)} - 0.1\chi.
\end{equation}
with $k=4$. The initial conditions in both the training and test sets are random draws from the curl of the `filtered velocity field' from JAX-CFD \cite{ml_accelerated_cfd}. 

The details of the standard solver are given by the `MUSCL' scheme described in \cref{sec:energyconservation}. The only difference between the machine learned solver and the standard solver is that the machine learned solver learns a correction to the flux at cell boundaries. This correction term is given by a simple CNN, which uses periodic convolutions and outputs 2 channels per cell. Each channel represents the flux at either the right cell boundary or the top cell boundary. The neural network uses kernels of size 5, has 64 channels per layer and 6 total layers.

We generate 100 evenly spaced data points from 100 simulations. For generation of training data, we run the `exact' high-resolution simulation until $t=20$ before sampling training data every $t=0.1$ units of time. Our loss function is again given by the mean squared error between the downsampled high-resolution `exact' time-derivative and the time-derivative from the machine learned solver. We train with a batch size of 100, use the Adam optimizer, and train for 1000 epochs over the training set. We use a learning rate of $10^{-4}$.

\Cref{fig:corr_2d_euler} is computed by computing the correlation between the downsampled high-resolution $128\times 128$ `exact' solution and the low-resolution solution(s). We run the `exact' solver until $t=20$ and use the solution at $t=20$ as the initial condition.

\noindent \textbf{Details for \cref{sec:verification_compressible_euler}}: Our goal is to solve the 1D compressible Euler's equations \cref{eq:compressibleeuler} and to compare the accuracy of a standard solver, a machine learned solver, and the same machine learned solver with the error-correcting invariant-preserving algorithm \cref{eq:compressible_euler_positivity_transformation,eq:compressible_euler_entropy_transformation} applied. We train solvers at different values of $N$, where $N$ is the number of spatial grid cells. We train solvers in both periodic domains and in domains with dirichlet boundary conditions. 

Our domain is $x \in [0, L]$ with $L=1$. The initial conditions in both the training and test sets are random draws from a relatively simple distribution. This distribution has $\rho_0 = \textnormal{max}\{\rho_{\textnormal{min}}, \rho_{s}\}$, $v_0 = v_{s}$, and $p_0 = \textnormal{max}\{p_{\textnormal{min}}, p_{s}\}$ where $\rho_s$, $v_s$, and $p_s$ are draws from (different) sine waves with random amplitudes and phases. The amplitude $A$ is drawn from a uniform distribution $A \in [0, 1]$. We set $\rho_{\textnormal{min}} = 0.75$ and $p_\textnormal{min} = 0.5$. 

We generate 20 snapshots from 10,000 different runs of a high-resolution `exact' simulation, for a total of 200,000 data points in the training set. Each run takes a snapshot at $t=0$ and then another snapshop every $t=0.01$ units of time until $t=0.2$. The exact simulation has $N=256$ grid cells. We train each solver with a batch size of 64 times the upsampling factor. Each solver is trained for 100,000 training iterations, a learning rate of $10^{-4}$, and uses the ADAM optimizer. 

The machine learned solver learns a correction to the fluxes of the standard solver at each cell boundary; the standard solver is once again given by the MUSCL scheme with reconstruction in characteristic variables.  The flux correction is given by the output of a simple CNN with 5 layers of 32 channels. Each layer pads with either periodic or edge padding. Each layer uses a kernel size of 5, except for the last layer which uses a kernel size of 4 due to symmetry considerations. 

Once again, our loss function is given by the mean squared error between the downsampled high-resolution `exact' time-derivative and the time-derivative from the machine learned solver.

\Cref{fig:1d_euler_verification} is computed by computing the MSE averaged over time and space and averaged over 50 different initializations in the test set.

\bibliography{cas-refs}
\bibliographystyle{elsarticle-num}





\end{document}